\documentclass[a4paper,11pt]{amsart}
\def\JOURNAL{1}  %

\usepackage{amsmath}
\usepackage{amsfonts}
\usepackage{amssymb}
\usepackage{listings}
\usepackage{mathtools}
\usepackage{hyperref}
\usepackage{amsthm}
\usepackage{verbatim}
\usepackage{appendix}
\usepackage[
]{todonotes}
\usepackage{url}

\if\JOURNAL1
\usepackage[a4paper, total={6in, 8in}]{geometry}
\fi

\if\JOURNAL2
\usepackage{natbib}
\usepackage{geometry}
\usepackage{hyperref}
\fi
\usepackage[T1]{fontenc}
\usepackage[latin1]{inputenc}
\usepackage{algorithm}
\usepackage{algorithmic}
\usepackage{array}
\usepackage{multirow}

\usepackage{subfigure}
\newcommand{\Rd}{{\mathrm{d}}}
\newcommand{\rset}{\mathbb{R}}

\newcommand{\TR}{\mathrm{Tr}\,} %

\newcommand{\xbar}{\bar{\mathbf{x}}}

\def\XXint#1#2#3{{\setbox0=\hbox{$#1{#2#3}{\int}$ }
\vcenter{\hbox{$#2#3$ }}\kern-.6\wd0}}

\numberwithin{figure}{section}
\numberwithin{table}{section}
\numberwithin{equation}{section}
\newtheorem{theorem}{Theorem}[section]

\newtheorem{preremark}[theorem]{Remark}
\if\JOURNAL1
\title[Path integral  molecular dynamics]{
Path integral molecular dynamics approximations of 
 quantum canonical observables}

\author[X. Huang]{Xin Huang}
\address{Institutionen f\"or Matematik, Kungl. Tekniska H\"ogskolan, 100 44 Stockholm, Sweden}
\email{xinhuang@kth.se}

\author[P. Plech\'a\v{c}] {Petr Plech\'a\v{c}}
\address{Department of Mathematical Sciences, University of Delaware, Newark, DE 19716, USA}
\email{plechac@udel.edu}

\author[M. Sandberg]{Mattias Sandberg}
\address{Institutionen f\"or Matematik, Kungl. Tekniska H\"ogskolan, 100 44 Stockholm, Sweden}
\email{msandb@kth.se}

\author[A. Szepessy]{Anders Szepessy}
\address{Institutionen f\"or Matematik, Kungl. Tekniska H\"ogskolan, 100 44 Stockholm, Sweden}
\email{szepessy@kth.se}
\fi

\begin{document}
\if\JOURNAL1
\maketitle

\begin{abstract}
Mean-field molecular dynamics based on path integrals is used to approximate canonical quantum observables for particle systems consisting of nuclei and electrons. A computational bottleneck is the sampling
from the Gibbs density of the electron operator, 
which due to the fermion sign problem 
has a computational complexity that scales exponentially with the number of electrons. 
In this work we construct an algorithm that approximates the mean-field Hamiltonian by path integrals for  fermions.  The algorithm is based on the determinant of a matrix with components based on Brownian bridges connecting  permuted  electron coordinates. The computational work for $n$ electrons is $\mathcal O(n^3)$, which reduces the computational complexity associated with the fermion sign problem. We 
analyze a bias resulting from this approximation and provide a computational error indicator. 
It remains to rigorously explain the surprisingly high accuracy. 
\end{abstract}
\tableofcontents

\fi
\if\JOURNAL2
\begin{frontmatter}

\title{%
Path integral molecular dynamics approximations of 
 quantum canonical observables}%

\author[1]{Xin Huang}
\ead{xinhuang@kth.se}

\author[2]{Petr Plech\'a\v{c}}
\ead{plechac@udel.edu}

\author[1]{Mattias Sandberg}
\ead{msandb@kth.se}

\author[1]{Anders Szepessy}
\ead{szepessy@kth.se}

\affiliation[1]{organization={Institutionen f\"or Matematik, Kungl. Tekniska H\"ogskolan},%
            city={Stockholm},
            postcode={100 44}, 
            country={Sweden}}

\affiliation[2]{organization={Department of Mathematical Sciences, University of Delaware},%
            city={Newark},
            postcode={19716},
            state={DE},
            country={USA}}

\begin{abstract}
Mean-field molecular dynamics based on path integrals is used to approximate canonical quantum observables for particle systems consisting of nuclei and electrons. A computational bottleneck is the sampling
from the Gibbs density of the electron operator, 
which due to the fermion sign problem 
has a computational complexity that scales exponentially with the number of electrons. 
In this work we construct an algorithm that approximates the mean-field Hamiltonian by path integrals for  fermions.  The algorithm is based on the determinant of a matrix with components based on Brownian bridges connecting  permuted  electron coordinates. The computational work for $n$ electrons is $\mathcal O(n^3)$, which reduces the computational complexity associated with the fermion sign problem. We 
analyze a bias resulting from this approximation and provide a computational error indicator. 
It remains to rigorously explain the surprisingly high accuracy. 
\end{abstract}

\begin{keyword}
ab initio molecular dynamics \sep canonical ensemble
\sep Gibbs distribution \sep path integral \sep fermion sign problem.
\MSC[2020] 35Q40 \sep 81S40 \sep 82C10 \sep 82M31.

\end{keyword}

\end{frontmatter}

\fi

\section{Background to approximations of quantum canonical observables}
Path integral methods are used to approximate observables
in the canonical quantum ensemble of particle systems consisting 
of nuclei and electrons. Specific implementations include, for instance,
ring polymer molecular dynamics, centroid molecular dynamics,
and related mean-field molecular dynamics, see \cite{tuckerman_compare} and \cite{HPSS}. 
A computational challenge arises especially for systems of fermions from the so-called (fermion) sign problem, 
causing the computational work of Monte Carlo approximations to grow exponentially with the number of fermion 
particles for a fixed expected approximation error \cite{dornheim_1, lyabartsev}. 
The aim of this work is to provide an alternative path integral formulation that reduces the sign problem,
where the Monte Carlo method averages the determinant of a matrix with components based on
Brownian bridge path integrals. 
The formulation uses the Feynman-Kac representation of the Gibbs partition function for the electron operator, 
with permuted Brownian bridge end coordinates due to 
indistinguishable particles. 
In this setting the partition function is represented by a sum over permutations of a rank four tensor. 
Sums of permutations of tensors are in general computationally feasible only in low dimension, since such sums 
often are NP-hard, see \cite{np_tensor}, %
with the exception of computing determinants of matrices. To reduce the tensor problem to determinants 
we use that the Gibbs potential energy for a nuclei-electron system can be written as a product based 
on path integrals for the energy with respect to each particle. %
In other words, we use that the potential energy is based on Coulomb pair interactions and its symmetry with respect to 
permutations of particles.
Our formulation has no bias in the case of a separable potential,
i.e., when it can be written as a sum of potentials that each depends on one particle only. %
Numerical tests show that the approximation error is relatively 
small also for non-separable potentials based on Coulomb electron 
repulsion, since the marginal distribution of initial and end particle positions turns out to be Gaussians.

We study the specific setting in which the electron kinetic energy operator is the Laplacian, 
which provides the Brownian bridge paths. Hence we consider the Hamiltonian
\[
\widehat{H}= -\frac{1}{2M}\Delta_{\mathbf X} - \frac{1}{2} \Delta_{\mathbf x} + V(\mathbf x,\mathbf X)=:-\frac{1}{2M}\Delta_{\mathbf X}+ H_e
\]
with  the Coulomb potential $V:\rset^{3n}\times\rset^{3N}\to\rset$,
the nuclei coordinates $\mathbf X=(\mathbf X_1,\ldots, \mathbf X_N)\in \rset^{3N}$,
and the electron coordinates $\mathbf x=(\mathbf x_1,\ldots\mathbf x_n)\in \rset^{3n}$, 
where $M$ is the nuclei-electron mass ratio.
We use the Hartree atomic units where the reduced Planck constant $\hbar=1$, the electron charge $e=1$, 
the Bohr radius $a_0=1$ and the electron mass $m_e=1$. 
The small semiclassical parameter $1/M$ is,
for example, in the case of a proton-electron system $1/M = m_e/m_p \approx 1/1836$.
The objective is to approximate canonical quantum correlation observables based 
on the Hamiltonian $\widehat{H}$ and the inverse temperature $\beta>0$
\[
\frac{\TR(\widehat A_t \widehat B_0 e^{-\beta\widehat{H}})}{\TR(e^{-\beta\widehat{H}})}
\]
and its symmetrized  form
\[
\mathfrak T_{{\rm qm}}:=\frac{\TR\big((\widehat A_t \widehat B_0+  \widehat B_0 \widehat A_t)e^{-\beta\widehat{H}}\big)}{2\TR(e^{-\beta\widehat{H}})}
\]
for quantum observables 
$\widehat A_t=e^{{\rm i} t \sqrt M \widehat{H}} \widehat A_0 e^{-{\rm i} t \sqrt M \widehat{H}}$ 
and $\widehat B_0$, at times $t$ and $0$, using  the trace, $\TR$, of a quantum operator over the nuclei and electron 
degrees of freedom. 
Such correlation observables are used, for instance, to determine the diffusion constant, reaction 
rates and other transport coefficients, \cite{Zwanzig_book}.

Our more precise aim is to obtain computable classical molecular dynamics approximations of  quantum observables.
For this purpose we need
the Weyl symbol $A:\rset^{6N}\to\mathbb C$ related to the quantum operator $\widehat A$, defined by
\begin{equation}\label{Weyl:definition}
\begin{split}
  \widehat{A} \phi(\mathbf X)=&
  \int_{\rset^{3N}} (\frac{\sqrt M}{2\pi})^{3N}\int_{\rset^{3N}} e^{{\rm i} M^{1/2}(\mathbf X-\mathbf Y)\cdot \mathbf P} 
   A\big(\frac{1}{2}(\mathbf X+\mathbf Y),\mathbf P\big)  \Rd \mathbf P\, \phi(\mathbf Y) \Rd \mathbf Y \\
\end{split}
\end{equation}
for Schwartz functions $\phi:\rset^{3N}\to\mathbb C$, see \cite{zworski}.
We assume that the Weyl symbols $A_0:\rset^{3N}\times\rset^{3N}\to \mathbb C$ and 
$B_0:\rset^{3N}\times\rset^{3N}\to \mathbb C$ for the initial quantum observables 
are independent of the electron coordinates.
We will use that the canonical quantum observables can be approximated by mean-field molecular dynamics
\begin{equation}\label{md_mf}
\mathfrak{T}_{\mathrm{md}}(t):=\frac{\int_{\rset^{6N}} A_0\big(\mathbf Z_t(\mathbf Z_0)\big) B_0(\mathbf Z_0) 
\TR( e^{-\beta H(\mathbf Z_0)})\Rd \mathbf Z_0}{\int_{\rset^{6N}}  \TR( e^{-\beta H(\mathbf Z_0)})\Rd \mathbf Z_0}
\end{equation}
where the Hamiltonian is given by %
\[
H=\frac{|\mathbf P|^2}{2} -\frac{1}{2}\Delta_{\mathbf x} + V(\mathbf x,\mathbf X)= \frac{|\mathbf P|^2}{2}+ H_e(\mathbf x, \mathbf X)
\]
and the trace, $\TR$, of a Weyl symbol is the trace with respect to the electron degrees of freedom, 
that is  the trace on an appropriate antisymmetric subset of $L^2(\rset^{3n})$.
The phase-space nuclei coordinates $\mathbf Z_t:=(\mathbf X_t,\mathbf P_t)$ solve the
mean-field Hamiltonian system
\begin{equation}\label{h}
\begin{split}
\dot{\mathbf X}_t =&\;\;\;\;\nabla_{\mathbf P} h(\mathbf X_t,\mathbf P_t)\,,\\
\dot{\mathbf P}_t =&-\nabla_{\mathbf X} h(\mathbf X_t,\mathbf P_t)\,,\\
\end{split}
\end{equation}
with the initial data $\mathbf Z_0:=(\mathbf X_0,\mathbf P_0)\in\rset^{3N}\times\rset^{3N}$ 
of nuclei positions and momenta, where the
mean-field Hamiltonian $h:\rset^{3N}\times \rset^{3N}\to\rset$ 
is defined by
\begin{equation}\label{mf-def}
h(\mathbf Z):=\frac{\TR (H(\mathbf Z)e^{-\beta H(\mathbf Z)})}{\TR (e^{-\beta H(\mathbf Z)})}\,.
\end{equation}
 Assuming that the electron problem for a given nuclei position, 
 \[
 H_e(\cdot,\mathbf X)\psi_i(\cdot,\mathbf X)=\lambda_i(\mathbf X)\psi_i(\cdot,\mathbf X)\,,
 \]  
 has the eigenvalues $\lambda_i(\mathbf X)\in \rset$, and eigenfunctions $\psi_i(\cdot,\mathbf X)\in L^2(\rset^{3n})$, $i=1,2,3,\ldots$,  then
\[
h(\mathbf X,\mathbf P)
=\frac{|\mathbf P|^2}{2} + \frac{\sum_{i=1}^\infty \lambda_i(\mathbf X)
e^{-\beta\lambda_i(\mathbf X)}
}{\sum_{i=1}^\infty e^{-\beta\lambda_i(\mathbf X)}}=:\frac{|\mathbf P|^2}{2} +\lambda_*(\mathbf X)\,.
\]

The work \cite{HPSS} derives the  approximation error 
 \begin{equation}\label{trace_est}
 |\mathfrak{T}_{\mathrm{qm}}(t)-\mathfrak{T}_{\mathrm{md}}(t)|=\mathcal O (M^{-1}+t\epsilon_1^2+t^2\epsilon_2^2)\,,
 \end{equation}
 where %
\[
\begin{split}
\epsilon_1^2 &= \frac{\|\TR\big((H-h)^2 e^{-\beta H}\big)\|_{L^1(\rset^{6N})}}{
\|\TR\big(e^{-\beta H}\big)\|_{L^1(\rset^{6N})}}
= \frac{\|\sum_{i=1}^\infty (\lambda_i-\lambda_*)^2
e^{-\beta\lambda_i}\|_{L^1(\rset^{3N})}
}{\|\sum_{i=1}^\infty e^{-\beta\lambda_i}\|_{L^1(\rset^{3N})}}\,,\\
\epsilon_2^2 &= 
\frac{\|\sum_{i=1}^\infty |\nabla(\lambda_i-\lambda_*)|^2
e^{-\beta\lambda_i}\|_{L^1(\rset^{3N})}
}{\|\sum_{i=1}^\infty e^{-\beta\lambda_i}\|_{L^1(\rset^{3N})}}\,,
\end{split}
\]
for the case that the electron part, $-\frac{1}{2}\Delta_{\mathbf x} + V(\mathbf x,\mathbf X)$, is approximated by
a finite dimensional matrix.
The mean-field Hamiltonian has the representation
\begin{equation}\label{h_simple}
\begin{split}
h(\mathbf X,\mathbf P)&=\frac{|\mathbf P|^2}{2}+\frac{\TR\Big( \big(-\frac{1}{2}\Delta_{\mathbf x}+V(\mathbf x,\mathbf X)\big)e^{-\beta(-\frac{1}{2}\Delta_{\mathbf x} + V(\mathbf x,\mathbf X))}\Big)}{\TR\big(e^{-\beta(-\frac{1}{2}\Delta_{\mathbf x} + V(\mathbf x,\mathbf X))}\big)}\\
&=\frac{|\mathbf P|^2}{2} + \frac{\TR (H_ee^{-\beta H_e})}{\TR (e^{-\beta H_e})}
= \frac{|\mathbf P|^2}{2} -\partial_\beta \log\big(\TR (e^{-\beta H_e})\big)\,.
\end{split}\end{equation}
To apply the molecular dynamics method \eqref{md_mf} requires
approximate sampling from the canonical Gibbs measure
\[
\mathbf X\mapsto\TR (e^{-\beta H_e(\cdot,\mathbf X)})/\int_{\rset^{3N}}\TR (e^{-\beta H_e(\cdot, \mathbf X_0)}){\mathrm{d}}\mathbf X_0:\rset^{3N}\to\rset\,.
\]
In particular, approximation of the mean-field $h:\rset^{6N}\to\rset$  is needed.
The purpose of this work is to formulate such approximations based on Monte Carlo sampled path integrals for fermion systems. 
The mean-field formulation \eqref{md_mf} is closely related to
the ring polymer molecular dynamics and centroid molecular dynamics methods, 
see \cite{ring-poly,tuckerman_compare,marx_hutter}.

The main result in this work is the construction 
and numerical tests of 
a proposed new computational method that approximates path integrals 
for fermions which is based on Brownian bridge paths and estimation of a suitable 
determinant in Section~\ref{sec_det}.
We present the necessary background on path integrals in Section~\ref{sec_path_integrals}. 
In Section~\ref{Motiv_Det}  we describe a perturbation analysis which, in combination with numerical experiments in Section~\ref{sec_experiment}, provides a computational error indicator to roughly
estimate the accuracy of the developed determinant formulation.
The numerical implementation and computational experiments for test cases with varying number of particles and inverse temperature are reported in Section~\ref{sec_experiment}.  
Finally, in Section~\ref{sec_conclusion} we comment  on a surprisingly high accuracy for which it remains 
to find a  rigorous theoretical answer.

\section{Path integrals}\label{sec_path_integrals}
First we present a brief background on path integral methods for the mean-field \eqref{h_simple},
to be used in Section~\ref{sec_det}
for the new representation of the partition function based on a determinant and Brownian bridge coordinate paths.

 We denote $\mathbf W:[0,\beta]\to\rset^{dn}$
 the standard Wiener process with independent components and define a Wiener path starting at $\mathbf{x}_0$ to be
$\mathbf x_t=\mathbf x_0+\mathbf W_t$.

\subsection{Distinguishable particles}
If electrons are treated as distinguishable particles
the Feynman-Kac  formulation  of path integrals, see
\cite[Theorem 7.6]{karatzas}  yields 
the expected value representation
\begin{equation}\label{FK1}
e^{-\beta H_e}\phi(\mathbf x_0)= \mathbb E[e^{-\int_0^\beta V(\mathbf x_t,\mathbf X){\mathrm{d}}t}\phi(\mathbf x_\beta)]\,.
\end{equation}
Then the $L^2(\rset^{dn})$ trace of the Gibbs density operator %
becomes
\begin{equation}\label{partition_d}
\TR(e^{-\beta H_e(\cdot,\mathbf X)})= \int_{\rset^{dn}} \mathbb E[e^{-\int_0^\beta V(\mathbf x_t,\mathbf X){\mathrm{d}}t}\delta(\mathbf x_\beta-\mathbf x_0)]{\mathrm{d}}\mathbf x_0\,.
\end{equation}

To approximate the mean-field Hamiltonian we compute
 \[
 \begin{split}
 \frac{\TR(H_e e^{-\beta H_e})}{\TR( e^{-\beta H_e})}
 &=-\partial_\beta \log \TR( e^{-\beta H_e})\\
& =-\partial_\beta \log \big(\frac{1}{(2\pi\beta)^{3n/2}}
\int_{\mathbb R^{3n}} \mathbb E[e^{-\int_0^\beta V(\mathbf x_t,\mathbf X){\mathrm{d}}t}
\, | \, \mathbf x_\beta=\mathbf x_0]{\mathrm{d}}\mathbf x_0\big)
 \end{split}
 \]
 where the conditional expectation is on Brownian bridge paths $\mathbf{x}_t$, $\mathbf{x}_\beta = \mathbf{x}_0$. 
 Hence, the integral $\int_0^\beta V(\mathbf x_t,\mathbf X){\mathrm{d}}t$ and the Brownian bridge $\mathbf x_t$ need to be discretized
 and generated. 
 One alternative is to directly generate discrete approximations of
 Brownian bridges, which is the focus in this work. 
 Another more common method, cf. \cite{ceperley3,dornheim_1}, is  
 to discretize the path integral as follows. Let $\mathbf x(j) := \mathbf x(j \beta/J)$ be the values of the paths at the time steps $t_j:= j\beta/J=:j\Delta t$. 
The normal distribution of the Brownian increments in the left hand side of \eqref{partition_d} implies that the Euler approximation representation,
 for the case with distinguishable particles
 \begin{equation}\label{Trotter_exp}
 \begin{split}
&\TR(e^{-\beta H_e(\cdot,\mathbf X)})\\
 &=\lim_{\Delta t\to 0+}\int_{\mathbb R^{dn}} 
\mathbb E[e^{-\sum_{j=0}^{J-1} V(\mathbf x(j),\mathbf X){\Delta t}}
\delta\big(\mathbf x(J) -\mathbf x(0)\big)]{\mathrm{d}}\mathbf x_0\\
&=\lim_{\Delta t\to 0+}\int_{\mathbb R^{dn}}\ldots \int_{\mathbb R^{dn}} \frac{e^{-\sum_{j=0}^{J-1} \big(V(\mathbf x(j),\mathbf X)\Delta t 
+\frac{|\mathbf x({j+1})-\mathbf x(j)|^2}{2\Delta t}\big)}}{(2\pi\Delta t)^{3nJ/2}}
\delta\big(\mathbf x(J) -\mathbf x(0)\big)
\prod_{j-0}^J \mathrm{d}\mathbf{x}(j)\\
%{\mathrm{d}}\mathbf x(0){\mathrm{d}}\mathbf x(1) \ldots {\mathrm{d}}\mathbf x(J)\\
%
%
%
&=\lim_{\Delta t\to 0+}\int_{\mathbb R^{dnJ}} \frac{e^{-\sum_{j=0}^{J-1} \big(V(\mathbf x(j),\mathbf X)\Delta t 
+\frac{|\mathbf x({j+1})-\mathbf x(j)|^2}{2\Delta t}\big)}}{(2\pi\Delta t)^{3nJ/2}}
{\mathrm{d}}\mathbf x(0){\mathrm{d}}\mathbf x(1) \ldots {\mathrm{d}}\mathbf x(J-1)\,,\\
\end{split}
\end{equation}
becomes the Trotter operator splitting, which is $\mathcal O((\Delta t)^2)$ accurate,  
see \cite{tretyakov_dumas}, due to the boundary condition $\mathbf x(J)=\mathbf x(0)$ in the right hand side. This integral in $\rset^{dnJ}$ is then approximated by the Metropolis or Langevin method applied to the potential 
\[
V_J(\mathbf x):=\sum_{j=0}^{J-1} \Big(V\big(\mathbf x(j),\mathbf X\big)\Delta t 
+\frac{|\mathbf x({j+1})-\mathbf x(j)|^2}{2\Delta t}\Big)
\] 
forming a classical ``ring polymer'' with  ``bead'' $\mathbf x(j)$, see
\cite[Chapter 10.2]{feynman-hibbs}.

\subsection{Indistinguishable particles}
In the case of $n$ indistinguishable  electrons, with the same spin,  we have instead the partition function
\begin{equation}\label{trace_1}
\TR(e^{-\beta H_e})= \int_{\rset^{dn}}\sum_{\sigma\in\mathbf S} \frac{{\rm sgn}(\sigma)}{n!}\mathbb E[e^{-\int_0^\beta V(\mathbf x_t,\mathbf X){\mathrm{d}}t}\delta(\mathbf x_\beta-\mathbf x_0^\sigma)]{\mathrm{d}}\mathbf x_0\,,
\end{equation} 
with the notation $\mathbf x^\sigma_t:=\big( \mathbf x_{\sigma_1}(t),\ldots,  \mathbf x_{\sigma_n}(t)\big)$ for 
$\mathbf x_t=\big(\mathbf x_1(t),\ldots, \mathbf x_n(t)\big)$ under any permutation $\sigma=(\sigma_1,\ldots,\sigma_n)$ of $(1,2,3,\ldots,n)$ and $t\in[0,\beta]$. The sum $\sum_{\sigma\in\mathbf S}$ is  over the set of all permutations $\mathbf S$ and $\mathbf x_k(t)\in\rset^d$ are the coordinates of electron $k$.
This expression of the partition function uses that the potential is invariant with respect to permutation of electron coordinates, i.e
$V(\mathbf x^\sigma,\mathbf X)=V(\mathbf x,\mathbf X)$ for all permutations $\sigma$. 
The setting with different spin is presented 
\if\JOURNAL1
in  Appendix~\ref{appendix_spin}.
\fi
\if\JOURNAL2
in  \ref{appendix_spin}.
\fi

The partition function \eqref{trace_1} can be derived by \eqref{FK1} using that the wave functions for the indistinguishable electrons have the antisymmetric representation 
$\phi(\mathbf x)=\sum_{\sigma\in\mathbf S}{\rm sgn}(\sigma)\tilde\phi(\mathbf x^\sigma)/n!$, for $\tilde\phi\in L^2(\rset^{dn})$, 
and consequently
\[e^{-\beta H_e}\phi(\mathbf x_0)=\mathbb E[\sum_{\sigma\in\mathbf S} \frac{{\rm sgn}(\sigma)}{n!} e^{-\int_0^\beta V(\mathbf x_t,\mathbf X){\mathrm{d}}t}\tilde\phi(\mathbf x_\beta^\sigma)]\,,\]
which implies \eqref{trace_1}.
In the case of bosons the wave functions are symmetric and the trace is replaced by 
\begin{equation}\label{boson}
\TR(e^{-\beta H_e})= \int_{\rset^{dn}}\sum_{\sigma\in\mathbf S} \frac{1}{n!}\mathbb E[e^{-\int_0^\beta V(\mathbf x_t,\mathbf X){\mathrm{d}}t}\delta(\mathbf x_\beta-\mathbf x_0^\sigma)]{\mathrm{d}}\mathbf x_0\,.
\end{equation} 

The partition function in \eqref{trace_1} can also be represented by Brownian bridges as follows.
A Brownian bridge process defined by
\begin{equation}\label{Def_Brownian_bridge}
\mathbf B_t:= \mathbf W_t - \frac{t}{\beta}\mathbf W_\beta\,,\quad 0\le t\le\beta\,,
\end{equation}
satisfies $\mathbf B_0=\mathbf B_\beta=0$
and is independent of $\mathbf W_\beta$, since the two Gaussians are uncorrelated
\[
\mathbb E[\mathbf B_t \mathbf W_\beta^*]
=\mathbb E[(\mathbf W_t - \frac{t}{\beta}\mathbf W_\beta) \mathbf W_\beta^*]
=0\,.
\]
The Wiener path $\mathbf x:[0,\beta]\to \rset^{3n}$ from $\mathbf x_0$ to
$\mathbf x_\beta=\mathbf x_0^\sigma$ in \eqref{trace_1}, can be represented by the Brownian bridge
\[
\begin{split}
\mathbf x_t&=
\mathbf x_0+\mathbf W_t
=\mathbf x_0+\frac{t}{\beta}\mathbf W_\beta +\mathbf B_t
=(1-\frac{t}{\beta})\mathbf x_0+\frac{t}{\beta} \mathbf x^\sigma_0 + \mathbf B_t\,.\\
\end{split}
\] 
The independence of $\mathbf B_t=\big(\mathbf B_1(t),\ldots,\mathbf B_n(t)\big)$ and $\mathbf W_\beta$ implies that
the partition function based on the Brownian bridge becomes
\begin{equation}\label{trace_bb}
\begin{split}
\TR(e^{-\beta H_e}) &= \int_{\rset^{3n}}\sum_{\sigma\in\mathbf S} \frac{{\rm sgn}(\sigma)}{n!}\mathbb E[e^{-\int_0^\beta V({\mathbf x}_t,\mathbf X){\mathrm{d}}t}\delta({\mathbf x}_\beta-\mathbf x_0^\sigma)]{\mathrm{d}}\mathbf x_0\\
&= \int_{\rset^{3n}}\sum_{\sigma\in\mathbf S} \frac{{\rm sgn}(\sigma)}{n!}
\mathbb E[e^{-\int_0^\beta V(\mathbf x_t,\mathbf X){\mathrm{d}}t}\ |\ \mathbf x_\beta=\mathbf x^\sigma_0]
\mathbb E[\delta(\mathbf W_\beta +\mathbf x_0-\mathbf x_0^\sigma)]{\mathrm{d}}\mathbf x_0\\
&= \int_{\rset^{3n}}\sum_{\sigma\in\mathbf S} \frac{{\rm sgn}(\sigma)}{n!}\mathbb E[e^{-
\int_0^\beta V(\mathbf x_t,\mathbf X){\mathrm{d}}t} \ |\ \mathbf x_\beta=\mathbf x^\sigma_0]\frac{e^{-|\mathbf x_0-\mathbf x_0^\sigma|^2/(2\beta)}}{(2\pi\beta)^{3n/2}}{\mathrm{d}}\mathbf x_0\,.\\
\end{split}
\end{equation}

The corresponding ring polymer method for fermions following \eqref{Trotter_exp} becomes
\begin{equation}\label{ring}
\begin{split}
&\TR(e^{-\beta H_e})\\
&=\lim_{\Delta t\to 0+}\sum_{\sigma\in\mathbf S} \frac{{\rm sgn}(\sigma)}{n!(2\pi\Delta t)^{3nJ/2}}
\int_{\mathbb R^{3n(J+1)}} e^{-V_J(\mathbf x)}\delta\big(\mathbf x(J)-\mathbf x^\sigma(0)\big)
{\mathrm{d}}\mathbf x(0){\mathrm{d}}\mathbf x(1) \ldots {\mathrm{d}}\mathbf x(J)\,.\\
\end{split}
\end{equation}

Monte Carlo approximations of the partition functions for bosons and fermions
are obtained by sampling both the permutations and the paths, 
based on the Wiener process.
In the case of low temperature, that is  $\beta\gg 1$, the partition function is dominated by the ground state eigenvalue,
which typically is substantially lower for bosons than for fermions. Consequently,
the value of the trace is lower for fermions than for bosons, 
while the variance in the Monte Carlo methods for the two remains similar 
unless the cancellations
for fermions can be removed by rewriting the fermions partition function as a  
sum with few cancelling terms. The fermion sign problem
refers to a much  larger computational work to approximate the fermion 
partition function by Monte Carlo samples of permutations and paths, 
compared to the boson partition function, in the case the sign cancellation 
is not explicitly included, cf. \cite{ceperley2, dornheim_1,holtzmann}.
The computationally costly sampling of all permutations has been improved by using the worm algorithm \cite{worm}.

The main new mathematical idea here is to write the fermion partition function \eqref{trace_bb} 
as an integral of a determinant
based on Brownian bridge paths. We thereby reduce the sign problem regarding the complexity with respect 
to the number of particles. The sign problem related to the increasing computational 
work with large $\beta$, however, remains.

Fermion determinants have been formulated before, both for quantum lattice models using hopping between
discrete lattice points in a second quantization setting, see the determinant quantum Monte Carlo method \cite{BSS,Li_Yao}, and for 
path integrals in the first quantization setting, based on the ring polymer approximations of \eqref{ring}, 
see \cite{takahashi, lyabartsev} and \eqref{Taka}. Our formulation also uses the path integral in the first quantization formulation \eqref{trace_bb} but
not the ring polymer approximation, instead we use a different determinant representation based on explicit generation of Brownian bridge paths from the Wiener process.
The first quantization formulation here, with Brownian bridge particle paths in $\rset^{3n}$, has the two advantages: (1) simple statistical independent sampling of Wiener processes as compared to Metropolis or Langevin sampling of ring polymers in high dimension,
and (2) the computational work, for $n$ particles requiring  $n$ classical paths, is proportional to $\mathcal O(n^3)$.
The drawback is that in our formulation 
two body interactions cannot be represented 
exactly. 
Nonetheless, we provide a computable error indicator based on perturbation analysis.

 For fermion particles in dimension one, $d=1$, the work \cite{ceperley3,dumas} shows that reflected Wiener processes can be used  to avoid the fermion sign problem in \eqref{trace_1}. Thereby the computational complexity becomes roughly
 $\mathcal O(\epsilon^{-3})$, as for distinguishable particles,
as a function of the expected approximation  error $\epsilon$.
However, for fermion particles in higher dimension, $d>1$, 
it seems that  the literature only provides computational methods to approximate \eqref{trace_1} 
with a computational complexity that grows exponentially with the number of particles $n$, 
cf. \cite{dornheim_1}, 
although there are different approximation alternatives, e.g.,
\cite{ceperley1,ceperley3,lyabartsev,parinello}, that improve the computational work compared 
to a straightforward Monte Carlo approximation of \eqref{trace_1}.

Simulations of bosons in \eqref{boson} do not experience the sign problem, nonetheless, they are also hard. 
For instance, the related problem to determine a matrix permanent is NP-hard, see \cite{np_tensor}. 

The fixed node method \cite{anderson, ceperley4, cances}
avoids the fermion sign problem for electron ground states, related to the formulation in $d=1$, 
but requires the nodal set of the wave function, which itself is challenging to determine and consequently 
causes computational bias that is difficult to estimate.

Our study relies on the electron operator, $H_e=-\frac{1}{2}\Delta_{\mathbf x} + V$,
and uses specific properties of  both 
the Laplacian  and the Coulomb interaction potentials for electrons and nuclei.

\section{The fermion partition functions based on determinants}\label{sec_det}

Takahashi and Imada \cite{takahashi} applied the fermion antisymmetrization to the density matrix approximation for each time step, 
namely to $e^{-(\frac{|\mathbf x(j+1)-\mathbf x(j)|^2}{2\Delta t} +V(\mathbf x(j),\mathbf X)\Delta t)}$ in \eqref{Trotter_exp},
 to obtain a formulation of the partition function based on products of determinants as follows
\begin{equation}\label{Taka}
\begin{split}
&\TR(e^{-\beta H_e})\\
&=\lim_{\Delta t\to 0+}\int_{\rset^{3nJ}} \prod_{j=0}^{J-1}\Big(
\sum_{\sigma\in\mathbf S} \frac{{\rm sgn}(\sigma)}{n!(2\pi\Delta t)^{3n/2}}
e^{-(\frac{|\mathbf x^\sigma(j+1)-\mathbf x(j)|^2}{2\Delta t} +V(\mathbf x(j),\mathbf X)\Delta t)}\Big)
{\mathrm{d}}\mathbf x(0)\ldots{\mathrm{d}}\mathbf x(J-1)\\
&=\lim_{\Delta t\to 0+}\int_{\rset^{3nJ}} 
\prod_{j=0}^{J-1}\Big( \frac{e^{-V(\mathbf x(j),\mathbf X)\Delta t}}{n!(2\pi\Delta t)^{3n/2}}{\rm det}\mathcal R\big(\mathbf x(j),\mathbf x(j+1)\big)\Big)
{\mathrm{d}}\mathbf x(0)\ldots{\mathrm{d}}\mathbf x(J-1)\,,\\
\end{split}
\end{equation}
where $\mathbf x(J)=\mathbf x(0)$ and the $n\times n$ matrix $\mathcal R$ has the components
\[
\mathcal R_{\ell k}\big(\mathbf x(j),\mathbf x(j+1)\big):=e^{|\mathbf x_{k}(j+1)-\mathbf x_\ell(j)|^2/(2\Delta t)}\,.
\]

A determinant can be computed with $\mathcal O(n^3)$  work. Therefore the formulation avoids the $\mathcal O(n!)$
computational complexity to sum over all permutations. However since the determinants have varying sign, for particles in dimension two and higher,
the required Monte Carlo sampler suffers from the fermion sign problem, which leads to a large statistical variance  \cite{takahashi,lyabartsev}.

In this section we present a different determinant formulation which is based on \eqref{trace_bb} where the Monte Carlo method directly 
samples independent
Brownian bridge paths, with the aim to reduce the fermion sign problem. 

The formulation of the partition function
\[
\TR(e^{-\beta H_e})
= \int_{\rset^{3n}}\sum_{\sigma\in\mathbf S} \frac{{\rm sgn}(\sigma)}{n!}\mathbb E[e^{-
\int_0^\beta V(\mathbf x_t,\mathbf X){\mathrm{d}}t} \ |\ \mathbf x_\beta=\mathbf x^\sigma_0]\frac{e^{-|\mathbf x_0-\mathbf x_0^\sigma|^2/(2\beta)}}{(2\pi\beta)^{3n/2}}{\mathrm{d}}\mathbf x_0
\]
using the Brownian bridge $\mathbf x_t=\big(1-\frac{t}{\beta}\big)\mathbf{x}_0+\frac{t}{\beta}\mathbf x_0^\sigma+\mathbf B_t$
can be represented with a determinant as follows.
Assume first that the potential is particle wise separable, so that
\begin{equation}\label{V_separated}
V(\mathbf x,\mathbf X) = \sum_{k=1}^n \tilde V_k(\mathbf x_k,\mathbf X),
\end{equation}
then
\[
H_e=-\frac{1}{2}\Delta_{\mathbf x} + V = \sum_{k=1}^n \big(-\frac{1}{2}\Delta_{\mathbf x_k} + \tilde V_k(\mathbf x_k,\mathbf X)\big)
\]
is separable. 
Define the $n\times n$ matrix 
\[
\mathcal W_{k\ell}\big(\mathbf x\big) := 
e^{-|\mathbf x_k(0)-\mathbf x_\ell(0)|^2/(2\beta)}
e^{-\int_0^\beta \tilde V_k(\mathbf B_k(t) + (1-\frac{t}{\beta}) \mathbf x_k(0)+ \frac{t}{\beta} \mathbf x_\ell(0),\mathbf X){\mathrm{d}}t}\,.
\]
Its determinant  yields the path integral
\[
\begin{split}
&{\mathrm{det}}\big(\mathcal W(\mathbf x)\big) =\sum_{\sigma\in\mathcal S} {\rm sgn}(\sigma) 
\mathcal W_{1\sigma_1}\mathcal W_{2\sigma_2}\ldots \mathcal W_{n\sigma_n}\\
&=\sum_{\sigma\in\mathcal S}{\rm sgn}(\sigma)e^{-\sum_{k=1}^n|\mathbf x_k(0)-\mathbf x_{\sigma_k}(0)|^2/(2\beta)}
e^{-\int_0^\beta\sum_{k=1}^n \tilde V_k(\mathbf B_k(t) + (1-\frac{t}{\beta}) \mathbf x_k(0)+ \frac{t}{\beta} \mathbf x_{\sigma_k}(0),\mathbf X){\mathrm{d}}t}\\
&=\sum_{\sigma\in\mathcal S}{\rm sgn}(\sigma)e^{-|\mathbf x_0-\mathbf x_0^\sigma|^2/(2\beta)}
e^{-\int_0^\beta  V(\mathbf B(t) + (1-\frac{t}{\beta}) \mathbf x(0)+ \frac{t}{\beta} \mathbf x^{\sigma}(0),\mathbf X) {\mathrm{d}}t}
\,,\\
\end{split}
\]
and we obtain 
the representation %
\begin{equation}\label{trace_det1}
\TR(e^{-\beta H_e}) = \int_{\rset^{3n}} \mathbb E[\frac{{\rm det}\big(\mathcal W(\mathbf x)\big) }{n!(2\pi\beta)^{3n/2}}]
{\mathrm{d}}\mathbf x_0\,.
\end{equation}

The path integral can be approximated by the 
trapezoidal method as follows.
Make a partition $t_m=m\Delta t$, for $m=0,\ldots M$,  with $\Delta t= \beta/M$
and let
\[
\begin{split}
\mathcal{\overline W}_{k\ell}\big(\mathbf x\big) &:= 
e^{-|\mathbf x_k(0)-\mathbf x_\ell(0)|^2/(2\beta)}
e^{-\sum_{m=1}^{M-1} \tilde V_k(\mathbf B_k(t_m) + (1-\frac{t_m}{\beta}) \mathbf x_k(0)+ \frac{t_m}{\beta} \mathbf x_\ell(0))\Delta t}\times\\
&\quad\times
e^{-\tilde V(\mathbf B_k(0) + \mathbf x_k(0))\Delta t/2
-\tilde V_k(\mathbf B_k(0) + \mathbf x_\ell(0))\Delta t/2
}\,.
\end{split}
\]
Then we have the computable second order accurate approximation 
\begin{equation}\label{trace_det2}
\TR(e^{-\beta H_e}) = \int_{\rset^{3n}} 
\mathbb E[\frac{{\rm det}\big(\mathcal{\overline W}(\mathbf x)\big) }{n!(2\pi\beta)^{3n/2}}]
{\mathrm{d}}\mathbf x_0 +\mathcal O((\Delta t)^2)\,.
\end{equation}
Such second order accuracy of path integral approximations is proved in \cite{tretyakov_dumas} assuming that the functional, here
$\mathbf x\mapsto e^{-\int_0^\beta V(\mathbf x_t){\mathrm{d}}t}$, is six times Fr{\'e}chet differentiable.
The determinant for an $n\times n$ matrix can be determined with $\mathcal O(n^3)$ number of operations  using the $LU$-factorization in contrast to the $\mathcal O(n!)$ computational work to sum over all permutations.

The Hamiltonian for molecular systems has a potential based on Coulomb interactions of the electrons and nuclei,
see \cite{handbook,lebris,lieb},
\begin{equation}\label{v_def}
\begin{split}
H_e(\mathbf x, \mathbf X)   & = -\frac{1}{2} \Delta_{\mathbf x} + V(\mathbf x,\mathbf X)\,,\\
V(\mathbf x,\mathbf X) & =\sum_{k=1}^n\sum_{k<k'} \frac{1}{|\mathbf x_k-\mathbf x_{k'}|} 
                      + \sum_{j=1}^N\sum_{j< i} \frac{Z_j Z_i}{|\mathbf X_j-\mathbf X_i |}
                       -\sum_{k=1}^n\sum_{j=1}^N \frac{ Z_j}{|\mathbf X_j-\mathbf x_k |}\,, %
\end{split}
\end{equation}
where  $Z_i$ denotes the charge of the $i$th nucleus. 
The coordinates $\mathbf x_k\in\rset^3$ and $\mathbf X_j\in \rset^3$ are for electron $k$ and nucleus $j$.
In this case the terms depending on $\mathbf x$ can still be  written as a sum over fermion particles
 \[
 \sum_{k=1}^n (\sum_{k'\ne k}\frac{1}{2|\mathbf x_k-\mathbf x_{k'}|}
 -\sum_{j=1}^N \frac{Z_j}{|\mathbf x_k-\mathbf X_j|})
 =: \sum_{k=1}^nV_k(\mathbf x,\mathbf X)
 \]
to form the potential energy per particle, $V_k$, but the electron repulsion sum
is not separable particle wise
and the repulsion, $\sum_{k'\ne k}\frac{1}{2|\mathbf x_k-\mathbf x_{k'}|}$, is instead
a sum  including all row coordinates $\mathbf x_{k'}$. 

Let 
\[
\mathbf x_{k}^{\sigma_k}(t):=\mathbf B_k(t) + (1-\frac{t}{\beta}) \mathbf x_k(0)+ \frac{t}{\beta} \mathbf x_{\sigma_k}(0)\,,
\]
then we have
\begin{equation}\label{tensor}
\begin{split}
\TR(e^{-\beta H_e}) &= \int_{\rset^{3n}}
\mathbb E \big[ \sum_{\sigma\in\mathbf S} \frac{{\rm sgn}(\sigma)}{n!(2\pi\beta)^{3n/2}}
e^{-\sum_{k=1}^n|\mathbf x_k(0)-\mathbf x_{\sigma_k}(0)|^2/(2\beta)}\times\\
&\quad \times e^{-\frac{\beta}{2}\sum_{i=1}^N\sum_{j=1,j\ne i}^N Z_iZ_j/|\mathbf X_i-\mathbf X_j|}\times\\
&\quad \times e^{\sum_{k=1}^n \int_0^\beta \sum_{j=1}^N Z_j/ |\mathbf x_k^{\sigma_k}(t) -\mathbf X_j| {\mathrm{d}} t}\times\\
&\quad \times e^{-\frac{1}{2}\sum_{k=1}^n \int_0^\beta \sum_{k'=1, k'\ne k}^n 
|\mathbf x_k^{\sigma_k}(t) - \mathbf x_{k'}^{\sigma_{k'}}(t)|^{-1} {\mathrm{d}} t} \big]
{\mathrm{d}}\mathbf x(0)\,,
\end{split}
\end{equation}
which cannot be written as a determinant of an $n\times n$ matrix, for $n>2$, due to the electron repulsion. Instead it is a sum over permutations for a tensor. Sums of permutations of tensors are often computationally demanding, e.g. $NP$-hard see \cite{np_tensor}.
An exception is the computation of the determinant for an $n\times n$ matrix which can be determined with $\mathcal O(n^3)$ operations.
The next step is therefore to formulate approximations to the right hand side in \eqref{tensor}, avoiding to sum over all permutations.

In the special case with two particles, $n=2$, the electron repulsion in \eqref{tensor} can in fact be
written as a determinant, since 
\begin{equation}\label{sigma_12}
\sigma_{k'}=\left\{
\begin{array}{cc}
2 & \mbox{ if } \sigma_k=1\,,\\
1 & \mbox{ if } \sigma_k=2\,.\\
\end{array}
\right.
\end{equation}
A natural approximation of \eqref{tensor} is to replace 
 $1/|\mathbf x_k^{\sigma_k}-\mathbf x_{k'}^{\sigma_{k'}}|$ 
by $1/|\mathbf x_k^{\sigma_k}-\mathbf x_{k'}^{\nu_{k'}}|$
where
\[
\nu_{k'} =\left\{\begin{array}{ll}
k', & \mbox{ if } k'\ne \sigma_k\,,\\
k, & \mbox{ if } k'=\sigma_{k}\,,
\end{array}\right.
\]
which for $n=2$ becomes the exact form \eqref{sigma_12} and
for $n>2$ treats the remaining particles with coordinates $\mathbf x_{k'}$ as distinguishable as follows. 
Define 
the approximation
\begin{equation}\label{Z_n}
\begin{split}
 \mathcal Z_\nu &:=\int_{\rset^{3n}}
\mathbb E \big[ \sum_{\sigma\in\mathbf S} \frac{{\rm sgn}(\sigma)}{n!(2\pi\beta)^{3n/2}}
e^{-\sum_{k=1}^n|\mathbf x_k(0)-\mathbf x_{\sigma_k}(0)|^2/(2\beta)}\times\\
&\quad \times e^{-\frac{\beta}{2}\sum_{i=1}^N\sum_{j=1, j\ne i}^N Z_iZ_j/|\mathbf X_i-\mathbf X_j|}\times\\
&\quad \times e^{\sum_{k=1}^n \int_0^\beta \sum_{j=1}^N Z_j/ |\mathbf x_k^{\sigma_k}(t) -\mathbf X_j| {\mathrm{d}} t}\times\\
&\quad \times e^{-\frac{1}{2}\sum_{k=1}^n \int_0^\beta \sum_{j=1,j\ne k}^n |\mathbf x_k^{\sigma_k}(t) - \mathbf x_{j}^{\nu_j}(t)|^{-1}{\mathrm{d}} t} \big]
{\mathrm{d}}\mathbf x(0)\,,
\end{split}
\end{equation}
which can be formulated by a determinant: let
\[
\begin{split}
\mathcal W_{k\ell}(\mathbf x) &:=
 e^{-|\mathbf x_k(0)-\mathbf x_{\ell}(0)|^2/(2\beta)}
 e^{\int_0^\beta \sum_{j=1}^N Z_j/ |\mathbf x_k^{\ell}(t) -\mathbf X_j| {\mathrm{d}} t}
e^{-\frac{1}{2} \int_0^\beta \sum_{j=1,j\ne k}^n |\mathbf x_k^{\ell}(t) - \mathbf x_{j}^{\nu_j}(t)|^{-1}{\mathrm{d}} t}\,,       
\end{split}
\]
where
\[
\nu_{j} =\left\{\begin{array}{ll}
j, & \mbox{ if } j\ne \ell\,,\\
k, & \mbox{ if } j=\ell\,,
\end{array}\right.\,
\]
then we have
\begin{equation}\label{Z_det}
 \mathcal Z_\nu=  e^{-\frac{\beta}{2}\sum_{i=1}^N\sum_{j=1, j\ne i}^N Z_iZ_j/|\mathbf X_i-\mathbf X_j|}\int_{\rset^{3n}} 
\mathbb E \big[\frac{{\rm det}\big(\mathcal W(\mathbf x)\big)}{n!(2\pi\beta)^{3n/2}}\big]
{\mathrm{d}}\mathbf x(0)\,,
\end{equation}
and $\mathcal Z_\nu=\TR(e^{-\beta H_e})$ for $n=2$. We also note that
for $\sigma$ equal to
the identity permutation or a permutation with only one transposition it holds that
\[
|\mathbf x_k^{\sigma_k}(t)-\mathbf x_j^{\sigma_j}(t)|=|\mathbf x_k^{\sigma_k}(t)-\mathbf x_j^{\nu_j}(t)|\,.
\]
For small values of $\beta$ such permutations will dominate
in \eqref{Z_n} and \eqref{tensor}, due to the factors 
$e^{-|\mathbf x_k(0)-\mathbf x_{\sigma_k}(0)|^2/(2\beta)}$, and consequently \eqref{Z_n}
is a consistent approximation of \eqref{tensor} as $\beta\to 0+$\,.

\begin{comment}
\textcolor{blue}{
We end this section with a comment on the determinant formulation which uses $\det\mathcal{R}$ in \eqref{Taka} from \cite{takahashi}.
The approximation of $\TR(e^{-\beta H_e})$
in \eqref{Taka} can be expressed in our framework 
using samples of Brownian bridges 
\[
\TR(e^{-\beta H_e}) \approx \lim_{\Delta t\to 0+}C'\int_{\rset^{3n}} \mathbb E\Big[
\prod_{j=0}^{J-1}\Big( \frac{e^{-V(\mathbf x(j),\mathbf X)\Delta t}}{
e^{-|\mathbf x(j+1)-\mathbf x(j)|^2/(2\Delta t)}}{\rm det}\mathcal R\big(\mathbf x(j),\mathbf x(j+1)\big)\Big)\Big]
{\mathrm{d}}\mathbf x(0)\,,
\]
where $\mathbf x(j)=\mathbf{x}(0)+\mathbf{W}(j\Delta t)-\frac{j}{J}\mathbf{W}(J\Delta t)$ 
and $\mathbb{E}[\cdot]$ is the expectation
under the standard Wiener process $\mathbf{W}$. 
In order to obtain a Monte Carlo sampling of \eqref{Taka} 
we have to include the 
density for the
Wiener increments in the denominator which renders 
this formulation impractical. Occurrence of exponentially small numbers in the denominator
prevents convergence of the empirical mean in Monte Carlo sampling.
\end{comment}

\subsection{Motivation for the determinant approximation \eqref{Z_n}}\label{Motiv_Det}

 The aim of this section is to motivate that for a given particle $k$ the
 interaction \[e^{-\frac{1}{2}\int_0^\beta \sum_{j\ne k}|\mathbf x_k^{\sigma_k}(t)
 - \mathbf x_j^{\sigma_j}(t)|^{-1}{\mathrm{d}}t}\] in \eqref{tensor}, where 
 $\mathbf x_i^{\sigma_i}(t)=\mathbf B_i(t)+(1-\frac{t}{\beta})\mathbf x_i(0)+\frac{t}{\beta}\mathbf x_{\sigma_i}(0)$, can be related  to the interaction in the determinant formulation 
 \[e^{-\frac{1}{2}\int_0^\beta \sum_{j\ne k}|\mathbf x_k^{\sigma_k}(t)
 - \mathbf x_j^{\nu_j}(t)|^{-1}{\mathrm{d}}t}\,.\] In particular we study the marginal distribution of $\mathbf x_{\sigma_j}(0)$ for  Brownian bridges from $\mathbf x(0)$ to $\mathbf x^\sigma(0)$. The mean of this marginal distribution motivates the choice 
 $\nu_j$.  This section shows that the variances of the marginals lead to
the perturbed interaction 
 \begin{equation}\label{xi_1}
\mathbb E_\xi[ e^{-\frac{1}{2}\int_0^\beta \sum_{j\ne k}|\mathbf x_k^{\sigma_k}(t)
 - (\mathbf x_j^{\nu_j}(t)+\frac{t}{\beta}\sqrt{\frac{\beta}{c_*}}\,\xi_j)|^{-1}{\mathrm{d}}t}]\,,\end{equation}
 which we use for sensitivity analysis of  the determinant formulation. Here the expected value is with respect  to only $\xi_j$, for $j\ne k$, which are independent
 standard normal distributed random variables in $\rset^3$ and $c_*\approx 2$ is  a parameter related to the logarithm of
 the permutation cycle length.  In particular, the numerical sensitivity experiments in Table \ref{Table_compare_perturb_BB_with_W_matrix} show that the perturbed mean-field approximation based on \eqref{xi_1} has a relative error on the same order of accuracy as $\bar h_\nu$ and thereby provides a computational
 error indicator %
 to roughly estimate the accuracy of $\bar h_\nu$, without using a reference value.

 \subsubsection{Marginals of $\mathbf x_{\sigma_j}(0)$} 
 If $\sigma={\rm I}$ we have $\mathbf x_j^{\sigma_j}=\mathbf x_j^{\nu_j}$
 and for small $\beta$ the permutation $\sigma={\rm I}$ dominates %
  due to the factor $e^{-|\mathbf x(0)-\mathbf x^{\sigma}(0)|^2/(2\beta)}$ in \eqref{tensor}.
  The case $\sigma={\rm I}$ implies to have only one-cycle permutations. In the case of a two cycle for particle $j$, i.e. $\sigma_{\sigma_j}=j$, the weight involving $\mathbf x_{\sigma_j}(0)$ becomes
  \[
  e^{-|\mathbf x_j(0)-\mathbf x_{\sigma_j}(0)|^2/(2\beta)} 
  e^{-|\mathbf x_{\sigma_j}(0)-\mathbf x_{\sigma_{\sigma_j}}(0)|^2/(2\beta)}
  =e^{-|\mathbf x_j(0)-\mathbf x_{\sigma_j}(0)|^2/\beta}\,.
  \]
  For a three-cycle, $\mathbf x_{\sigma_{\sigma_{\sigma_j}}}(0)=\mathbf x_j(0)$,
  we have 
  \[\begin{split}
&e^{-|\mathbf x_j(0)-\mathbf x_{\sigma_j}(0)|^2/(2\beta)} 
  e^{-|\mathbf x_{\sigma_j}(0)-\mathbf x_{\sigma_{\sigma_j}}(0)|^2/(2\beta)}
e^{-|\mathbf x_{\sigma_{\sigma_j}}(0)-\mathbf x_{\sigma_{\sigma_{\sigma_j}}}(0)|^2/(2\beta)}\\
&=e^{-|\mathbf x_j(0)-\mathbf x_{\sigma_j}(0)|^2/(2\beta)}
e^{-|\mathbf x_{\sigma_{\sigma_j}}(0)-\frac{\mathbf x_j(0)+\mathbf x_{\sigma_j}(0)}{2\beta}|^2/\beta}
e^{-|\mathbf x_j(0)-\mathbf x_{\sigma_j}(0)|^2/(4\beta)}\,,
 \end{split} 
  \]
  using that for any $y,y_a$ and $y_b$ in $\rset^3$
  \[
  e^{-|y-y_a|^2/\gamma_a}e^{-|y-y_b|^2/\gamma_b}
  =e^{-|y-\frac{\gamma_by_a+\gamma_a y_b}{\gamma_a\gamma_b}|^2
  \frac{\gamma_a+\gamma_b}{\gamma_a\gamma_b}}
  e^{-|y_a-y_b|^2/(\gamma_a+\gamma_b)}\,.
  \]
  Integration with respect to $\mathbf x_{\sigma_{\sigma_j}(0)}$ yields the non normalized marginal distribution
  \[
  e^{-|\mathbf x_j(0)-\mathbf x_{\sigma_j}(0)|^2(\frac{1}{2\beta}+\frac{1}{4\beta})}
  =e^{-|\mathbf x_j(0)-\mathbf x_{\sigma_j}(0)|^2\frac{3}{4\beta}}\,.
  \]
  For a four-cycle with $\mathbf x_{\sigma_{\sigma_{\sigma_{\sigma_j}}}}(0)=\mathbf x_j(0)$ we obtain similarly the weight
  \[
  e^{-|\mathbf x_j(0)-\mathbf x_{\sigma_j}(0)|^2/(2\beta)} 
  e^{-|\mathbf x_{\sigma_j}(0)-\mathbf x_{\sigma_{\sigma_j}}(0)|^2/(2\beta)}
e^{-|\mathbf x_{\sigma_{\sigma_j}}(0)-\mathbf x_{\sigma_{\sigma_{\sigma_j}}}(0)|^2/(2\beta)}
e^{-|\mathbf x_{\sigma_{\sigma_{\sigma_j}}}(0)-\mathbf x_{\sigma_{\sigma_{\sigma_{\sigma_j}}}}(0)|^2/(2\beta)}
  \]
  which by integration with respect to $\mathbf x_{\sigma_{\sigma_j}}(0)$ and 
  $\mathbf x_{\sigma_{\sigma_{\sigma_j}}}(0)$ determines the marginal distribution proportional to
  \[
  e^{-|\mathbf x_j(0)-\mathbf x_{\sigma_j}(0)|^2(\frac{1}{2\beta}+\frac{1}{4\beta}+\frac{1}{4\beta+2\beta})}
 = e^{-|\mathbf x_j(0)-\mathbf x_{\sigma_j}(0)|^2\frac{11}{12\beta}}\,.
  \]
For an $m$-cycle, $m\ge 3$, we obtain analogously the marginal distribution proportional to
\begin{equation}\label{normal_x}
e^{-|\mathbf x_j(0)-\mathbf x_{\sigma_j}(0)|^2
\sum_{m'=1}^{m-1}\frac{1}{2m'\beta}}\,,
\end{equation}
where  
\begin{equation}\label{m_l}
\sum_{m'=1}^{m-1}\frac{1}{2m'\beta}=:\frac{c_{m}}{2\beta}\approx \frac{c_*}{2\beta}\,,
\end{equation}
and $c_{2}=2, c_3=3/2, c_4=11/6, c_5=25/12\,.$ The approximation
$c_*=2$  for small $\beta$ is further motivated in Section \ref{sub_section_perturbation}.

In the determinant formulation \eqref{Z_det}
\[
\mathcal W_{k\ell} = 
e^{-|\mathbf x_k(0)-\mathbf x_{\ell}(0)|^2/(2\beta)}
 e^{\int_0^\beta \sum_{j=1}^N Z_j/ |\mathbf x_k^{\ell}(t) -\mathbf X_j| {\mathrm{d}} t}
e^{-\frac{1}{2} \int_0^\beta \sum_{j=1,j\ne k}^n |\mathbf x_k^{\ell}(t) - \mathbf x_{j}^{\nu_j}(t)|^{-1}{\mathrm{d}} t}\,,
\]
we do not have access to the permutations. Therefore we approximate $\mathbf x_{\sigma_j}(0)$ in
the interaction $|\mathbf x_k^{\sigma_k}(t)
 - \mathbf x_j^{\sigma_j}(t)|^{-1}$
by the mean $\mathbf x_j(0)$ of the normal distributions
$e^{-|\mathbf x_j(0)-\mathbf x_{\sigma_j}(0)|^2\sum_{m'=1}^{m-1}\frac{1}{2m'\beta}}$,
with respect to the variable $\mathbf x_{\sigma_j}(0)$, related to $m$-cycles.
If $j=\ell$ and $\ell\ne k$ we cannot use $\nu_j=\ell$
since then $|\mathbf x_k^\ell(\beta)-\mathbf x_\ell^\ell(\beta)|=0$
and consequently we would obtain $\mathcal W_{k\ell}=0$  for all $\ell$ in the numerical approximations.
Instead we let $\nu_\ell=k$ in order for $\nu_j,\, j\ne k$, to all be different.

\subsubsection{The perturbation formulation \eqref{xi_1}} \label{sub_section_perturbation}
To motivate the  formulation \eqref{xi_1}, we introduce the notation
\[
\begin{split}
\TR(e^{-\beta H_e}) &=\sum_{\sigma\in\mathbf S} \frac{{\rm sgn}(\sigma)}{n!(2\pi\beta)^{3n/2}}
 \int_{\rset^{3n}} \mathbb E[e^{-
\int_0^\beta V(\mathbf x_t,\mathbf X){\mathrm{d}}t} \ |\ \mathbf x_\beta=\mathbf x^\sigma_0]
e^{-|\mathbf x_0-\mathbf x_0^\sigma|^2/(2\beta)}{\mathrm{d}}\mathbf x_0\\
&=:\sum_{\sigma\in\mathbf S} {\rm sgn}(\sigma)
\int_{\rset^{3n}} f_\sigma(\mathbf x_0)
e^{-|\mathbf x_0-\mathbf x_0^\sigma|^2/(2\beta)}
{\mathrm{d}}\mathbf x_0\,.\\
\end{split}
\]
Assume that for a permutation $\sigma$ the coordinate $\mathbf x_j(0)$ is in an $m$-cycle
consisting of the coordinates 
\[\big(\mathbf x_j(0),\mathbf x_{\sigma_j}(0),\ldots, \mathbf x_{\sigma_{\ldots\sigma_j}}(0)\big)
=:y_\sigma\in\rset^{3m}
\]
and let the particle coordinates that are not in this $m$-cycle be denoted by $y^\perp_\sigma\in\rset^{3(n-m)}$. Define 
\[
f_\sigma\big(\mathbf x_1(0),\mathbf x_2(0),\ldots,\mathbf x_n(0)\big)
=:  f_{\sigma,j}(y_\sigma,y^\perp_\sigma)\,.
\]
We have by \eqref{normal_x}
\[
\begin{split}
&\int_{\rset^{3n}} f_\sigma(\mathbf x_0)
e^{-|\mathbf x_0-\mathbf x_0^\sigma|^2/(2\beta)}{\mathrm{d}}\mathbf x_0\\
&=\int_{\rset^{3n}}f_{\sigma,j}(y_\sigma,y^\perp_\sigma) 
e^{-|\mathbf x_j(0)-\mathbf x_{\sigma_j}(0)|^2/(2\beta)}e^{-|\mathbf x_{\sigma_j}(0)-\mathbf x_{\sigma_{\sigma_j}}(0)|^2/(2\beta)}
e^{-|\mathbf x_{\sigma_{\sigma_j}}(0)-\mathbf x_{\sigma_{\sigma_{\sigma_j}}}(0)|^2/(2\beta)}\ldots\times\\
&\qquad\times e^{-|y^\perp-(y^{\perp})^\sigma|^2/(2\beta)}{\mathrm{d}}\mathbf x_j(0){\mathrm{d}}\mathbf x_{\sigma_j}(0)\ldots{\mathrm{d}}\mathbf x_{\sigma_{\ldots\sigma_j}}(0)
{\mathrm{d}}y^\perp_\sigma 
\\
&=C_{mn}\int_{\rset^{3(n+1-m)}}f_{\sigma,j}\Big(\mathbf x_j(0),\mathbf x_{\sigma_j}(0),
\bar{\mathbf x}_{\sigma_{\sigma_j}}(0),\bar{\mathbf x}_{\sigma_{\sigma_{\sigma_j}}}(0),
\ldots,\bar{\mathbf x}_{\sigma_{\ldots\sigma_j}}(0), y^\perp_\sigma\Big)\times\\ 
&\qquad\times e^{-|\mathbf x_j(0)-\mathbf x_{\sigma_j}(0)|^2\sum_{m'=1}^{m-1}{\frac{1}{2m'\beta}}}
e^{-|y^\perp-(y^{\perp})^\sigma|^2/(2\beta)}{\mathrm{d}}\mathbf x_j(0){\mathrm{d}}\mathbf x_{\sigma_j}(0)
{\mathrm{d}}y^\perp_\sigma\,,
\end{split}
\]
where $\big(\bar{\mathbf x}_{\sigma_{\sigma_j}}(0),\bar{\mathbf x}_{\sigma_{\sigma_{\sigma_j}}}(0),\ldots,\bar{\mathbf x}_{\sigma_{\ldots\sigma_j}}(0)\big)$ and the constant $C_{mn}$
are defined by the integration with respect to 
$\big(\mathbf x_{\sigma_{\sigma_j}}(0),\ldots ,\mathbf x_{\sigma_{\ldots\sigma_j}}(0)\big)\in\rset^{3(m-2)}\,,
$
using the mean value theorem. The mean value theorem also yields
\begin{equation}\label{x_bar}
\begin{split}
    &\int_{\rset^{3n}} f_\sigma(\mathbf x_0)
\frac{e^{-|\mathbf x_0-\mathbf x_0^\sigma|^2/(2\beta)}}{(2\pi\beta)^{3n/2}}{\mathrm{d}}\mathbf x_0\\
&=C'_{mn}\int_{\rset^{3(n+1-m)}}f_{\sigma,j}\Big(\mathbf x_j(0),
\bar{\mathbf x}_{\sigma_j}(0),
\bar{\mathbf x}'_{\sigma_{\sigma_j}}(0),\bar{\mathbf x}'_{\sigma_{\sigma_{\sigma_j}}}(0),
\ldots,\bar{\mathbf x}'_{\sigma_{\ldots\sigma_j}}(0),y^\perp_\sigma\Big)\times\\ 
&\qquad\times 
e^{-|y^\perp-(y^{\perp})^\sigma|^2/(2\beta)}{\mathrm{d}}\mathbf x_j(0)
{\mathrm{d}}y^\perp_\sigma\,,
\end{split}
\end{equation}
by integration with respect to the Gaussian measure $e^{-|\mathbf x_j(0)-\mathbf x_{\sigma_j}(0)|^2\sum_{m'=1}^{m-1}{\frac{1}{2m'\beta}}}{\mathrm{d}}\mathbf x_{\sigma_j}(0)$,
 for another constant $C'_{mn}$ and  new values $(\bar{\mathbf x}'_{\sigma_{\sigma_j}}(0),\bar{\mathbf x}'_{\sigma_{\sigma_{\sigma_j}}}(0)),
\ldots,\bar{\mathbf x}'_{\sigma_{\ldots\sigma_j}}(0))$ obtained from the mean value theorem and the dependence on $\mathbf x_{\sigma_j}(0)$. 

The work \cite{dornheim_c} shows that, provided the temperature is high enough, $\beta\lesssim 1$, the probability of cycle length $m$ with high probability decays exponentially 
for the ideal Fermi gas, i.e. for non interacting electrons, and for the uniform electron gas, i.e. for interacting electrons. On the other hand, in the case of low temperature,  $\beta= 8$, Figure 9 in \cite{dornheim_c} shows that the probability is almost constant with respect to the cycle length.
In our study focused on high temperature, we therefore approximate the sum in \eqref{normal_x}  and \eqref{m_l} as
\begin{equation*}\label{normal_x2}
e^{-|\mathbf x_j(0)-\mathbf x_{\sigma_j}(0)|^2
\sum_{m'=1}^{m-1}\frac{1}{2m'\beta}} 
\approx 
e^{-|\mathbf x_j(0)-\mathbf x_{\sigma_j}(0)|^2 \frac{c_*}{2\beta}
}\,,
\end{equation*}
where $c_*=2$.
The dependence on $\sigma$ in \eqref{normal_x} and \eqref{m_l}, through the cycle length $m$, is thereby avoided so that we can apply a determinant formulation.
The normal distribution of $\mathbf x_{\sigma_j}(0)$ obtained in \eqref{normal_x}  and \eqref{x_bar} 
provides a method to study perturbations of the determinant formulation \eqref{Z_det}, where the size of the perturbations are related to the unavailable 
interactions $|\mathbf x_k^{\sigma_k}(t)-\mathbf x_j^{\sigma_j}(t)|$ for the determinant formulation. 
Therefore we use in \eqref{xi_1} the $\sigma$ independent approximation
\begin{equation}\label{xi_approx}
\bar{\mathbf x}_{\sigma_j}(0)\approx \mathbf x_j(0)+\sqrt{\frac{\beta}{c_*}} \, \xi_j
\end{equation}
where $\xi_j\in\rset^3$ are independent standard normal random variables. 
A corresponding rigorous analysis of the perturbations would require to more precisely determine the approximation error in \eqref{xi_approx} using the independent samples $\xi_j$.
The perturbation approximation \eqref{xi_approx}   
motivates to use the following error indicator for the determinant formulation
\begin{equation}\label{Z-perturbed}
\begin{split}
&|\TR(e^{-\beta H_e}) -\mathcal Z_\nu| \\
&=C\Big|
\sum_{\sigma\in\mathcal S}{\rm sgn}(\sigma)\int_{\rset^{3n}}
\mathbb E\Big[\prod_{k=1}^n\Big( e^{\int_0^\beta \sum_{j=1}^N Z_j/ |\mathbf x_k^{\sigma_k}(t) -\mathbf X_j| {\mathrm{d}} t}\times\\
 &\qquad\times
e^{-\sum_{j\ne k}\frac{1}{2}\int_0^\beta 
|\mathbf x_k^{\sigma_k}(t)-\mathbf x_j^{j}(t) -(\mathbf x_j^{\sigma_j}(t)-
\mathbf x_j^{j}(t))|^{-1}{{\mathrm{d}}t}}
e^{-|\mathbf x_k(0)-\mathbf x_{\sigma_k}(0)|^2/(2\beta)}\Big)\Big]{\mathrm{d}}\mathbf x_0
\\ &\quad 
-\int_{\rset^{3n}} \mathbb E[{\rm det} \big(\mathcal W(\mathbf x)\big)]
{\mathrm{d}}\mathbf x_0\Big|
\\
&\approx C 
|\int_{\rset^{3n}} \mathbb E[{\rm det}\big( \widetilde{\mathcal W}(\mathbf x)\big)]
{\mathrm{d}}\mathbf x_0
-\int_{\rset^{3n}} \mathbb E[{\rm det} \big(\mathcal W(\mathbf x)\big)]
{\mathrm{d}}\mathbf x_0|
\end{split}
\end{equation}
where
\begin{equation}\label{W_bar}
\begin{split}
\widetilde{\mathcal W}_{k\ell}(\mathbf x_0) &=
 e^{\int_0^\beta \sum_{j=1}^N Z_j/ |\mathbf x_k^{\ell}(t) -\mathbf X_j| {\mathrm{d}} t}\times\\
 &\qquad\times
\mathbb E_\xi[e^{-\sum_{j\ne k}\frac{1}{2}\int_0^\beta 
{|\mathbf x_k^\ell(t)-(\mathbf x_j^{\nu_j}(t) 
+ \frac{t}{\beta}\sqrt{\frac{\beta}{c_*}}\xi_j)|^{-1}}{{\mathrm{d}}t}}]
e^{-|\mathbf x_k(0)-\mathbf x_\ell(0)|^2/(2\beta)}\,,\\
C&:=\frac{e^{-\frac{\beta}{2}\sum_{i=1}^N\sum_{j=1, j\ne i}^N Z_iZ_j/|\mathbf X_i-\mathbf X_j|}}{n!(2\pi\beta)^{3n/2}}\,.
\end{split}
\end{equation}
Here the expected value of ${\rm det}(\widetilde{\mathcal W})$ 
is with respect to $\mathbf B$ and the expeded value $\mathbb E_\xi$ is with respect to $\xi$ only.

\begin{preremark}[Bias in mean-field approximations]
In the case of separable potentials the determinant mean-field formulation \eqref{Z_det} and 
the Hartree-Fock mean-field method for electron ground state approximations both have no bias.
The determinant formulation also has no bias in the case of interacting 
distinguishable particles, which is not the case for the Hartree method. Consequently the determinant formulation is asymptotically accurate in the two important settings of separable fermion potentials and for interacting distinguishable particles.
\end{preremark}

\section{Numerical experiments}\label{sec_experiment}
In this section, we apply the path integral 
with Feynman-Kac  formulation \eqref{trace_det2} 
and \eqref{Z_det}, %
to approximate the quantum statistical partition function
and the mean-field energy for systems consisting of a few indistinguishable electrons. These systems operate at various inverse temperatures $\beta$ and involve two different potential functions denoted as $V(\mathbf{x},\,\mathbf{X})$, namely:
    \begin{subequations}
     \begin{align}
     \tag{V1}
      V(\mathbf{x},\mathbf{X}) &= \frac{1}{2}|\mathbf{x}|^2,\textup{ Section }\ref{Subsection:Fermi_gas} ,\label{Case_V1}\\
      \tag{V2}
      V(\mathbf{x},\mathbf{X}) &= \frac{1}{2}|\mathbf{x}|^2+\frac{1}{2}\sum_{k}\sum_{\ell\ne k}\frac{\lambda}{|x_k-x_{\ell}|},\textup{ Section }\ref{Subsection_Coulomb_repulsion_plus_HO}. \label{Case_V2}
     \end{align}
    \end{subequations}

    The derived partition function approximations facilitate the computation of corresponding mean-field energy values. To validate our approach, we compare our results for the two test cases, \ref{Case_V1} and \ref{Case_V2}, with three different benchmarks: exact solutions obtained for particle-wise separable potentials, reference data in \cite{dornheim_1}, and the formula \eqref{trace_bb} employing all the permutations. Furthermore, we implement the perturbation formulation in Section \ref{sub_section_perturbation} to calculate an error indicator for our numerical results. The MATLAB code for implementing the algorithm based on \eqref{trace_det2} and \eqref{Z_det} is openly available through our GitHub repository \cite{PIMC_github_repo}.
\subsection{Non-interacting Fermi gas in a confining harmonic trap} \label{Subsection:Fermi_gas}
\medskip
In $d$-dimensional space with a governing separable external potential, we study a system comprising $n$ indistinguishable non-interacting fermions, known as a Fermi gas. Our investigation involves a test case where the Fermi gas is confined by an external harmonic trap defined by $V_{\mathrm{HO}}(\mathbf{x})=\frac{1}{2}|\mathbf{x}|^2$. To determine the partition function's exact value at inverse temperature $\beta=1\,/\,(k_B T)$, we use a recursive formula, given by
\begin{equation}\label{trace_recursion}
    Z_n(\beta) = \frac{1}{n}\,\sum_{k=1}^{n}\,(-1)^{k-1}\,Z_1(k\beta)\,Z_{n-k}(\beta)\,,\mbox{ with }Z_0(\beta)=1\,,
\end{equation}
where 
\begin{equation}\label{one_particle_Z}
     Z_1(\beta)= \sum_{j=1}^\infty\,e^{-\beta\,\epsilon_j}\,,
\end{equation}
represents the single-particle partition function. Here, $\epsilon_j$ denotes the energy levels of each eigenstate corresponding to %
$V_{\mathrm{HO}}(\mathbf{x})$, see \cite{Borrmann_trace_recursion,Schmidt_trace_recursion}. We compute the exact partition function using \eqref{trace_recursion} as the reference value, which serves to test the convergence of the proposed path integral Monte Carlo method \eqref{trace_det2}.

Considering the simple case of a one-dimensional external harmonic oscillator potential, the single-particle energy eigenvalues in Hartree atomic units are expressed as $\eta_j = \frac{1}{2}+j$, where $j=0,1,2,\dots$. Using these eigenvalues, we can directly deduce the single-particle partition function $Z_1(\beta)$ for a three-dimensional harmonic oscillator trap:
\begin{equation}\label{single_particle_Z1_3d_HO}
 Z_1(\beta)=\sum_{j,k,\ell=0}^\infty e^{-\beta(\eta_j+\eta_k+\eta_{\ell})}=\Big( \sum_{j=0}^\infty e^{-\beta\eta_j} \Big)^3=\Big( \frac{\exp{(-\frac{\beta}{2})}}{1-\exp{(-\beta)}} \Big)^3\,. 
\end{equation}
By combining the recursive formula \eqref{trace_recursion} with \eqref{single_particle_Z1_3d_HO}, we can compute the exact partition function $Z_n(\beta)$. In contrast, using \eqref{Z_det}, we can approximate the partition function for the test case without nuclei by
\begin{equation}\label{partition_function_HO_approx}
    \Bar{Z}_n(\beta) := \frac{1}{M_x}\sum_{m=1}^{M_x} \frac{\mathrm{det}\big( \overline{\mathcal{W}}(\mathbf{x}(\cdot\,;m)) \big)}{n!\,(2\pi\beta)^{3n/2}\,p(\mathbf{x}(0\,;m))},
\end{equation}
which is the Monte Carlo approximation of the multi-dimensional integral 
\[\int_{\rset^{3n}} \mathbb E\Big[\frac{{\rm det}\big(\mathcal{\overline W}(\mathbf x)\big) }{n!(2\pi\beta)^{3n/2}}\Big] {\mathrm{d}}\mathbf x_0\]
based on $M_x$ samples of independent paths 
\begin{equation} \label{x_t_BB_path}
\mathbf{x}(t\,;i)=\mathbf{B}(t\,;i)+\mathbf{x}(0\,;i), \quad\textup{for }i=1,\dots,M_x,
\end{equation}
 where $\mathbf{x}(0\,;i)\in\mathbb{R}^{3n}$ is the $i$-th sample of the initial position, and $\mathbf{B}(t\,;i)\in\mathbb{R}^{3n}$ is an independent Brownian bridge process sample defined by \eqref{Def_Brownian_bridge}, which satisfies $\mathbf{B}(0)=\mathbf{B}(\beta)=0$. The samples of the initial positions $\{ \mathbf{x}(0\,;i) \}_{i=1}^{M_{x}}$ are independently drawn from a mixed normal distribution, with probability density function
\begin{equation}\label{p_x_density}
p(\mathbf{x})=\frac{1}{2}\Big( \frac{1}{\sigma_1\sqrt{2\pi}}\,e^{-\frac{|\mathbf{x}|^2}{2\sigma_1^2}} + \frac{1}{\sigma_2\sqrt{2\pi}}\, e^{-\frac{|\mathbf{x}|^2}{2\sigma_2^2}}\Big)\,,
\end{equation}
where $\sigma_1^2=\beta$ and $\sigma_2^2 = 1/\beta$ are parameters adjusting the variance of the sampling 
distribution for various inverse temperatures $\beta$.
The matrix element in the $k$-th row and ${\ell}$-th column of $\mathcal{\overline W}(\mathbf{x})$ is given by
\[
\begin{aligned}
\mathcal{\overline W}_{k\ell}\big(\mathbf x\big) &= 
e^{-|\mathbf x_k(0)-\mathbf x_\ell(0)|^2/(2\beta)}
e^{-\sum_{m=1}^{M-1} \tilde V_k(\mathbf B_k(t_m) + (1-\frac{t_m}{\beta}) \mathbf x_k(0)+ \frac{t_m}{\beta} \mathbf x_\ell(0))\Delta t}\times\\
&\quad\times
e^{-\big(\tilde V_k( \mathbf x_k(0))+\tilde V_k( \mathbf x_\ell(0))\big)\Delta t/2
}\,.
\end{aligned}
\]
Here we use the trapezoidal method to evaluate the integration in the imaginary time range $[0,\beta]$ numerically, with step size $\Delta t=\beta/M$, discrete time points $t_m = m\,\Delta t$ for $m=0,1,\dots, M$, and $\mathbf{B}_k(t)$ is the $k$-th component of the Brownian bridge process $\mathbf{B}(t)$ . Due to the absence of pairwise interactions between the fermions in this model, the potential function $V(\mathbf{x},\mathbf{X})$ depends only on the external harmonic oscillator term and is particle-wisely separable
\[ V(\mathbf{x},\mathbf{X})=\sum_{k=1}^n \Tilde{V}_k(\mathbf{x}_k),\ \textup{ with }\ \Tilde{V}_k(\mathbf{x}_k)= V_{\mathrm{HO}}(\mathbf{x}_k)=\frac{1}{2}|\mathbf{x}_k|^2\,. \]
Specifically, we use different time-step sizes $\Delta t$ to obtain the approximations $\Bar{Z}_n(\beta)$ for $n=6$ non-interacting fermions at varying $\beta$ values. The results of our numerical implementations are summarized in Table ~\ref{Table_approx_Z_HO}. 
\begin{table}[h]\label{Table_approx_Z_HO}
\begin{center}
 \begin{tabular}{ c c c c c c c }
\hline\hline\\ [-2.3ex]\multicolumn{6}{c}{\ }
\\[-2.3ex]
 $\beta$ & $M_{x}$ & $Z_n(\beta)$ &  $\Delta t$ & $\Bar{Z}_n(\beta)$ & Relative diff. & Relative CI  \\ \multicolumn{5}{c}{\ }\\[-2.3ex]
 \hline
\multicolumn{6}{c}{\ } \\ [-1.0ex]
\multirow{2}{*}{$1.0$} & \multirow{2}{*}{$2^{28}$} & \multirow{2}{*}{$1.6978 \times 10^{-4}$} & $0.025$ & $1.6983(8)\times 10^{-4}$& $2.96 \,\mathrm{e}(-4)$ & $4.54 \,\mathrm{e}(-4)$\\
 &  &  & $0.0125$ & $1.6980(8)\times 10^{-4}$ & $1.08 \,\mathrm{e}(-4)$ & $4.56 \,\mathrm{e}(-4)$\\
 [-1.0ex]\multicolumn{6}{c}{\ } \\
\multirow{2}{*}{$1.5$} & \multirow{2}{*}{$2^{28}$} & \multirow{2}{*}{$5.9174\times 10^{-9}$} & $0.025$ & $5.922(13)\times 10^{-9}$ & $8.45 \,\mathrm{e}(-4)$ & $2.14 \,\mathrm{e}(-3)$\\

 &  &  & $0.0125$ & $5.915(13)\times 10^{-9}$ & $3.43 \,\mathrm{e}(-4)$ & $2.14 \,\mathrm{e}(-3)$\\
[-1.0ex]\multicolumn{6}{c}{\ } \\
\multirow{2}{*}{$2.0$} & \multirow{2}{*}{$2^{28}$} & \multirow{2}{*}{$6.8663\times 10^{-13}$} & $0.025$ & $7.1(5)\times 10^{-13}$& $3.59 \,\mathrm{e}(-2)$ & $7.81\,\mathrm{e}(-2)$\\
 &  &  & $0.0125$ & $6.9(5)\times 10^{-13}$ & $7.25 \,\mathrm{e}(-3)$ & $7.65 \,\mathrm{e}(-2)$\\
[-1.25ex]
\multicolumn{5}{c}{\ } \\ \hline\hline
\end{tabular} 
\end{center}
\caption{  %
Case \ref{Case_V1}: Approximation of partition function $\Bar{Z}_n(\beta)$ obtained by \eqref{partition_function_HO_approx} for $n=6$ non-interacting fermions in dimension 3, compared with the exact value from \eqref{trace_recursion}.  The digit in the parenthesis of the $\Bar{Z}_n(\beta)$ values denotes the statistical error in our numerical result by the Monte Carlo integration method.
The relative difference between $\Bar{Z}_n(\beta)$ and $Z_n(\beta)$ is computed by $|\Bar{Z}_n(\beta)-Z_n(\beta)|/Z_n(\beta)$. The 95\% relative confidence interval of $\Bar{Z}_n(\beta)$ is recorded in the last column.
The total sample size of the independent Monte Carlo estimators for evaluating the multi-dimensional integral in $\mathbb{R}^{3n}$ is denoted by $M_{x}$. 
}
\end{table}
\smallskip
Furthermore, we also obtain estimations for the system mean-field energy $h(\beta)$ based on the derivatives of the approximated partition function $\Bar{Z}_n(\beta)$ 
with respect to the parameter $\beta$
\[ 
h(\beta)=-\partial_\beta \log\big(Z_n(\beta)\big)\simeq -\partial_\beta \log{\big(\Bar{Z}_n(\beta)\big)}=:\Bar{h}(\beta)\,. 
\]
We apply a rescaling technique on the sampled Brownian bridge paths to obtain the derivatives of the matrix element $\partial_\beta \mathcal{W}_{k\ell}$. %
This technique allows a formulation using integration based on standard Brownian bridge samples over the time range $[0,1]$, avoiding the need for evaluating the derivatives of the Brownian bridge processes with respect to $\beta$ over the time range $[0,\beta]$.
Further details regarding this rescaling technique are provided in the beginning of Section~\ref{Subsection_Coulomb_repulsion_plus_HO} 
\if\JOURNAL1
and in Appendix~\ref{Appendix_section_a}. 
\fi
\if\JOURNAL2
and ~\ref{Appendix_section_a}. 
\fi

To establish a reference value for the mean-field energy $h(\beta)$, we utilize a simple central difference formula to approximate the derivative:
\[
\partial_\beta \log\big(Z_n(\beta)\big)\simeq \frac{\log{\big(Z_n(\beta+\Delta\beta)\big)}-\log{\big(Z_n(\beta-\Delta\beta)\big)}}{2\Delta\beta}+\mathcal{O}\big( \Delta\beta\, ^2\big)\,.
\]
Here, $Z_n(\beta+\Delta\beta)$ and $Z_n(\beta-\Delta\beta)$ are computed using the exact recursive formula \eqref{trace_recursion}, and an inner loop gradually reduces the $\Delta\beta$ value to attain the desired accuracy for the reference value of $h(\beta)$.

\begin{figure}[h]
\centering

\subfigure[Error in $\Bar{Z}_n(\beta)$.]{
  \includegraphics[width=0.45\textwidth]{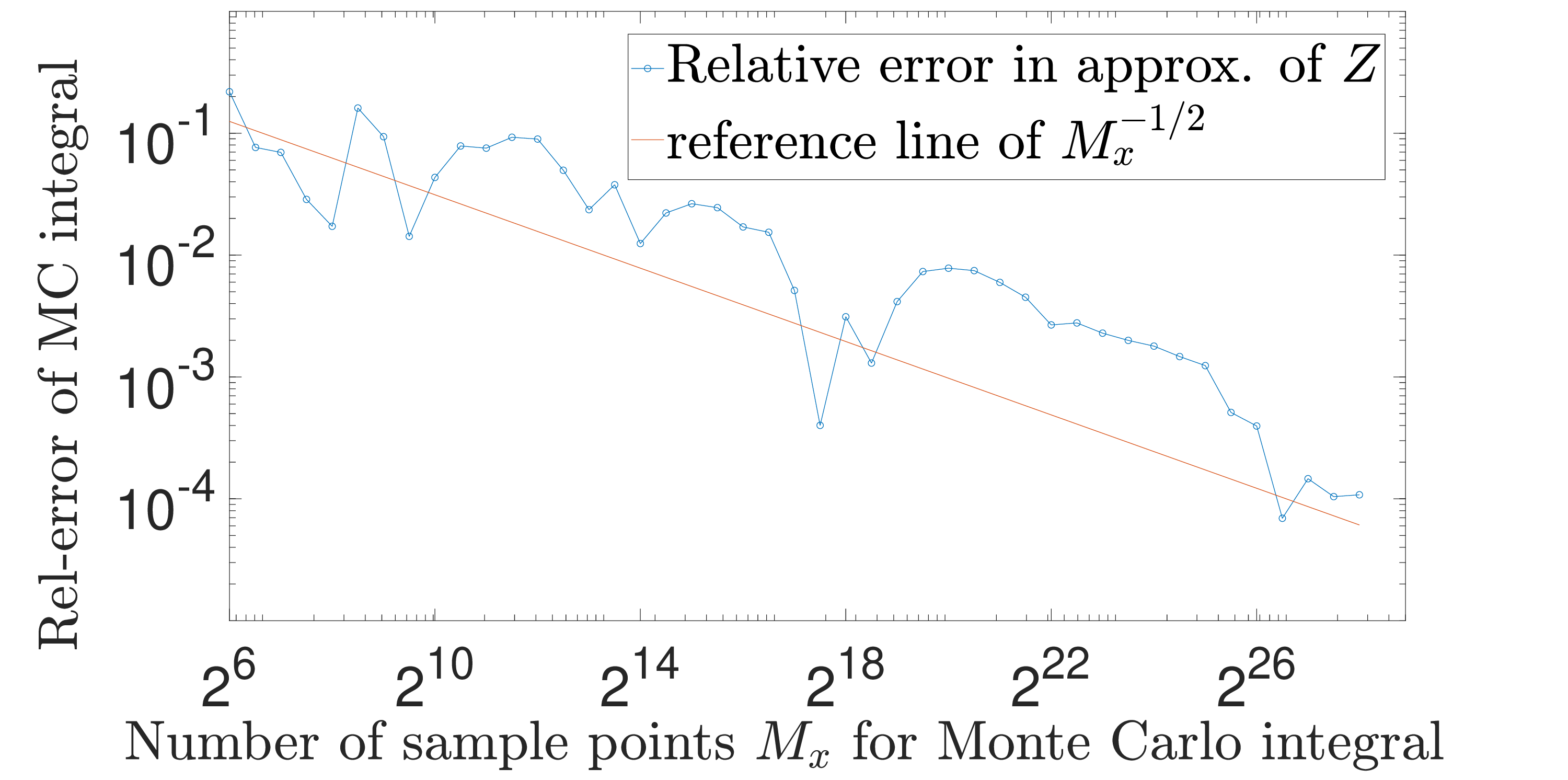}
  \label{fig:pure_HO_approx_Z}
}
\hspace*{-0.3cm}
\subfigure[Error in $\Bar{h}(\beta)$.]{
  \includegraphics[width=0.45\textwidth]{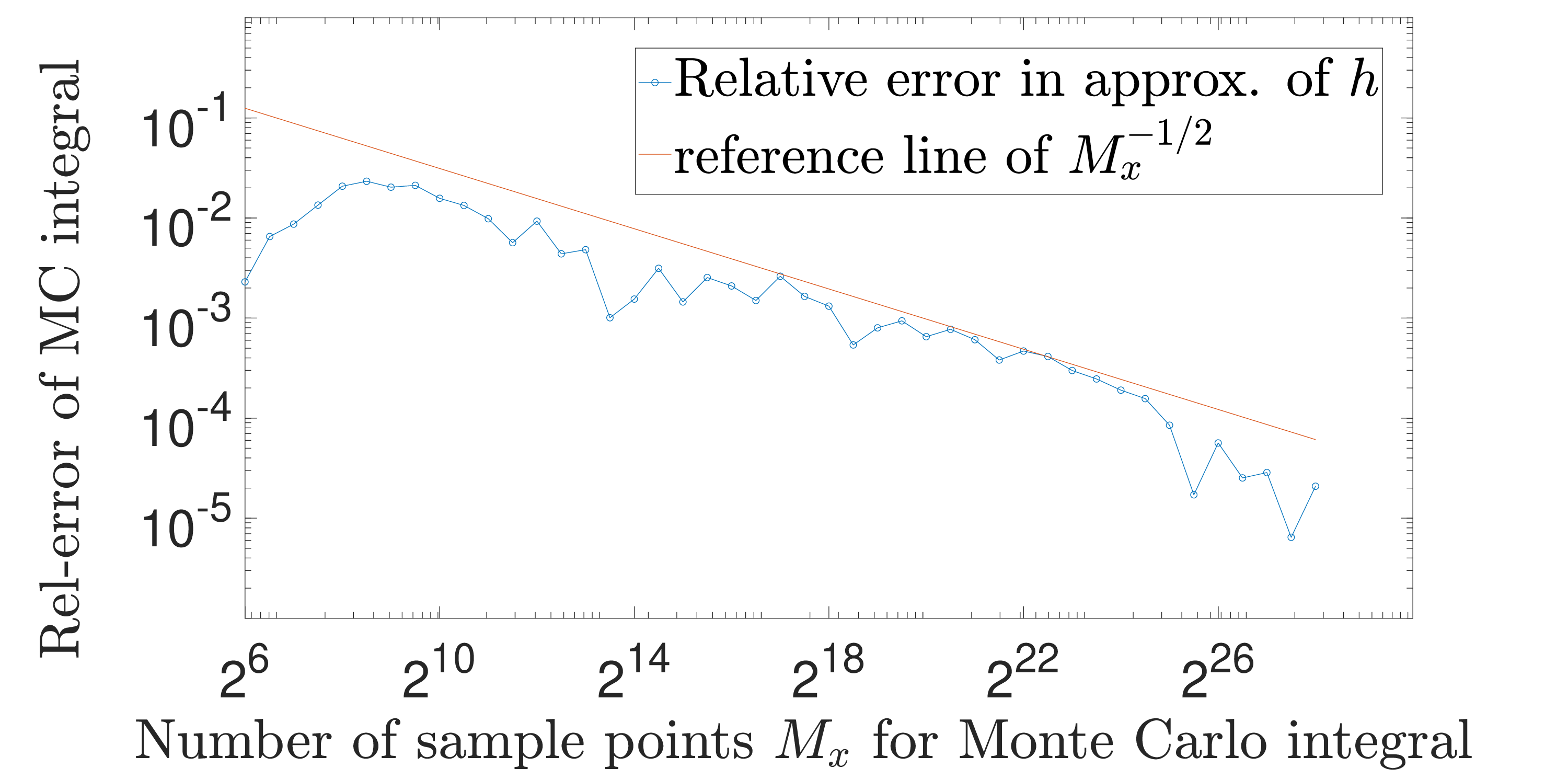}
  \label{fig:pure_HO_approx_h}
}

\caption{Case \ref{Case_V1}: The convergence rate of the errors in our Monte Carlo approximation of the partition function $\Bar{Z}_n(\beta)$ and the mean-field energy $\Bar{h}(\beta)$ for the fermi gas model under external harmonic oscillator potential, with increasing sample size $M_{x}$. For both subfigures, the parameters are chosen as inverse temperature $\beta = 1$, number of fermions $n=6$, dimensionality $d=3$, and the time step size $\Delta t=0.0125$}.
\label{Fig:rel_error_pure_HO}
\end{figure}

The corresponding results for our approximation of the system mean-field energy $\Bar{h}(\beta)$ are summarized in
Table~\ref{Table_approx_mean-field_HO}. 
In Figure~\ref{Fig:rel_error_pure_HO}, we plot the error of our Monte Carlo approximations to the partition function $\bar{Z}_n(\beta)$ and to the mean-field energy $\Bar{h}(\beta)$ with increasing sample size $M_{x}$ at inverse temperature $\beta = 1$.
\begin{table}[h]
\begin{center}
 \begin{tabular}{ c c c c c c c }
\hline\hline\\ [-2.3ex]\multicolumn{6}{c}{\ }
\\[-2.3ex]
 $\beta$ & $M_{x}$ & $h(\beta)$ &  $\Delta t$ & $\Bar{h}(\beta)$ & Relative diff. & Relative CI  \\ \multicolumn{5}{c}{\ }\\[-2.3ex]
 \hline
\multicolumn{6}{c}{\ } \\ [-1.0ex]
\multirow{2}{*}{$1.0$} & \multirow{2}{*}{$2^{28}$} & \multirow{2}{*}{$22.7799$} & $0.025$ & $22.778(10)$ & $7.81 \,\mathrm{e}(-5)$ & $4.44 \,\mathrm{e}(-4)$\\
 &  &  & $0.0125$ & $22.780(10)$ & $2.09 \,\mathrm{e}(-5)$ & $4.45 \,\mathrm{e}(-4)$\\
 [-1.0ex]\multicolumn{6}{c}{\ } \\
\multirow{2}{*}{$1.5$} & \multirow{2}{*}{$2^{28}$} & \multirow{2}{*}{$18.9572$} & $0.025$ & $18.95(3)$ & $2.14 \,\mathrm{e}(-4)$ & $1.66 \,\mathrm{e}(-3)$\\

 &  &  & $0.0125$ & $18.96(3)$ & $4.86 \,\mathrm{e}(-5)$ & $1.66 \,\mathrm{e}(-3)$\\
[-1.0ex]\multicolumn{6}{c}{\ } \\
\multirow{2}{*}{$2.0$} & \multirow{2}{*}{$2^{28}$} & \multirow{2}{*}{$17.4894$} & $0.025$ & $17.2(8)$ & $1.53 \,\mathrm{e}(-2)$ & $4.44 \,\mathrm{e}(-2)$\\
 &  &  & $0.0125$ & $17.4(8)$ & $5.88 \,\mathrm{e}(-3)$ & $4.54 \,\mathrm{e}(-2)$\\
[-1.25ex]
\multicolumn{5}{c}{\ } \\ \hline\hline
\end{tabular} 
\end{center}
\caption{  %
Case \ref{Case_V1}: Approximation of the mean-field energy $\Bar{h}(\beta)$ for $n=6$ non-interacting fermions in dimension 3, compared with the corresponding exact value $h(\beta)$. The digit in the parenthesis of the $\Bar{h}(\beta)$ values denotes the statistical error by the Monte Carlo integration method.
The relative difference between $\Bar{h}(\beta)$ and $h(\beta)$ is computed by $|\Bar{h}(\beta)-h(\beta)|/h(\beta)$, and the 95\% relative confidence interval of $\Bar{h}(\beta)$ is recorded in the last column.
The total sample size of the independent Monte Carlo estimators for evaluating the multi-dimensional integral in $\mathbb{R}^{3n}$ is denoted by $M_{x}$. 
}
\label{Table_approx_mean-field_HO}
\end{table}

From Tables~\ref{Table_approx_Z_HO} and \ref{Table_approx_mean-field_HO}, we observe that the statistical error of our Monte Carlo estimation of the high-dimensional integration in the Feynman-Kac formulation is dominating the time discretization error since the difference between using $\Delta t=0.025$ and $\Delta t=0.0125$ is smaller than the corresponding confidence interval based on $M_{x}=2^{28}$ samples. The statistical error of the Monte Carlo estimation grows also quickly as the inverse temperature $\beta$ increases, due to the more severe cancellation of permuted particle configurations at higher $\beta$ value.

Moreover, since no additional approximations are introduced for this example with the absence of pairwise interactions between particles, we have that $\Bar{Z}_n(\beta)$ and $\Bar{h}(\beta)$ will approach the exact partition function $Z_n(\beta)$ and mean-field energy $h(\beta)$, in the limit of $\Delta t\to 0$ and $M_{x}\to\infty$. This is also manifested by the error curves 
in Figure~\ref{Fig:rel_error_pure_HO}, where no biases are observed as the sample size $M_{x}$ increases. Also, the errors decrease with a rate approximately proportional to $M_x^{-\frac{1}{2}}$, confirming the well-known convergence rate of Monte Carlo approximation of the high-dimensional integral. As a comparison, we also calculate the partition function for a system consisting of distinguishable particles. The computationally faster approximation of the fermionic mean-field energy by distinguishable particle partition function gives an error of $15\%$ for $\beta = 1$ and $25\%$ for $\beta = 1.5$.

\subsection{Electrons with Coulomb interaction under a harmonic trap} \label{Subsection_Coulomb_repulsion_plus_HO}
\medskip
For a model including nuclei and involving pairwise Coulomb interactions, our approximation to the system partition function reads
\[\mathcal{Z}_\nu=e^{-\frac{\beta}{2}\sum_{i=1}^N\sum_{j=1,j\neq i}^N Z_iZ_j/|\mathbf{X}_i-\mathbf{X}_j|}\int_{\mathbb{R}^{dn}}\mathbb{E}\big[ \frac{\det{\big( \mathcal{W}(\beta,\mathbf{x})\big)}}{n!(2\pi\beta)^{dn/2}} \big]\,\mathrm{d}\mathbf{x}(0)\,,\]
as given in \eqref{Z_det}. Recalling the corresponding definition of the matrix element %
\[
\mathcal W_{k\ell}(\beta,\mathbf x) :=
 e^{-|\mathbf x_k(0)-\mathbf x_{\ell}(0)|^2/(2\beta)}
 e^{\int_0^\beta \sum_{j=1}^N Z_j/ |\mathbf x_k^{\ell}(t) -\mathbf X_j| {\mathrm{d}} t}
e^{-\frac{1}{2} \int_0^\beta \sum_{j=1,j\ne k}^n 1/|\mathbf x_k^{\ell}(t) - \mathbf x_{j}^{\nu_j}(t)|{\mathrm{d}} t}\,,       
\]
where 
\[
\nu_{j} =\left\{\begin{array}{ll}
j, & \mbox{ if } j\ne \ell\,,\\
k, & \mbox{ if } j=\ell\,,
\end{array}\right.\,
\]
we can take the derivative of $\mathcal{Z}_{\nu}$ with respect to the inverse temperature $\beta$, and obtain the corresponding approximation of the system mean-field energy $h$ by
\[
h:=\frac{\TR (H_e e^{-\beta H_e})}{\TR (e^{-\beta H_e})} \approx -\frac{\partial_\beta\, \mathcal{Z}_{\nu}(\beta)}{\mathcal{Z}_{\nu}(\beta)}
\,.
\]
To evaluate the derivative $\partial_\beta\, \mathcal{Z}_{\nu}(\beta)$, we need to have the derivative of the determinant value $\partial_\beta\det{\big( \mathcal{W}(\beta,\mathbf{x})\big)}$. %
Jacobi's formula yields
\begin{equation} \label{Jacobi_formula}
\frac{\partial}{\partial\beta}\mathrm{det}\big(\mathcal{W}(\beta,\mathbf{x})\big)= 
\mathrm{Tr}\Big( \mathrm{adj}\big(\mathcal{W}(\beta,\mathbf{x})\big)\, \frac{\partial \mathcal{W}(\beta,\mathbf{x}) }{\partial \beta} \Big)\,,
\end{equation}
where $\frac{\partial \mathcal{W}(\beta,\mathbf{x}) }{\partial \beta}$ is the element-wise derivative of matrix $\mathcal{W}(\beta,\mathbf{x})$ with respect to $\beta$, and $\mathrm{adj}\big(\mathcal{W}(\beta,\mathbf{x})\big)$ denotes the adjugate of the matrix $\mathcal{W}(\beta,\mathbf{x})$, see \cite{Stewart_def_adjugate}. The $(i,j)$-element of this adjugate matrix is defined by the transpose of the co-factor matrix of $\mathcal{W}(\beta,\mathbf{x})$, i.e.,
\[ \mathrm{adj}\big(\mathcal{W}(\beta,\mathbf{x})\big)_{ij} :=  (-1)^{i+j}\det{\big( M_{ji} \big)} , \]
with $M_{pq}$ denoting the submatrix of $\mathcal{W}(\beta,\mathbf{x})$ obtained by deleting it $p$-th row and $q$-th column.

We consider a system consisting of $n$ identical electrons, with pairwise Coulomb repulsive interactions under a confining harmonic oscillator potential. The matrix element $\mathcal{W}_{k\ell}(\beta,\mathbf{x})$ takes the form
\begin{equation}\label{W_mat_element_no_nuclei}
\mathcal W_{k\ell}(\beta,\mathbf x) =
 e^{-\frac{|\mathbf x_k(0)-\mathbf x_{\ell}(0)|^2}{2\beta}} 
e^{ -\int_0^\beta \Tilde{V}_{k\ell}(\mathbf{x}(t)){\mathrm{d}} t}\,,       
\end{equation}
where
\begin{equation}\label{V_tilde_nu_j}
\Tilde{V}_{k\ell}\big(\mathbf{x}(t)\big)=\frac{1}{2} |\mathbf{x}_k^{\ell}(t)|^2 + \sum_{j=1,j\ne k}^n \frac{\lambda}{2|\mathbf{x}_k^{\ell}(t) - \mathbf{x}_j^{\nu_j}(t)|}\,, 
\end{equation}
and the parameter $\lambda$ measures the relative strength of the Coulomb interactions compared to the harmonic oscillator potential. Throughout the numerical experiments  of this subsection, we choose the relative force parameter $\lambda=0.5$. Here $\mathbf{x}_k^{\ell}(t)$ denotes a Brownian bridge path with initial position $\mathbf{x}_k(0)$ at $t=0$ and ending position $\mathbf{x}_{\ell}(0)$ at $t=\beta$, defined by
\begin{equation}\label{x_k_ell(t)_path_def_1}
\mathbf x_{k}^{\ell}(t):=\mathbf B_k(t) + (1-\frac{t}{\beta}) \mathbf x_k(0)+ \frac{t}{\beta} \mathbf x_{\ell}(0)\,,\ t\in[0,\beta]\,,
\end{equation}
where $\mathbf{B}_k(t)$ denotes a Brownian bridge with $\mathbf{B}_k(0)=\mathbf{B}_k(\beta)=0$, as defined in \eqref{Def_Brownian_bridge}.

In order to determine $\frac{\partial \mathcal{W}(\beta,\mathbf{x})}{\partial \beta}$, it is useful to define the stochastic process $\overline{\mathbf{B}}:[0,1]\times\Omega\to\mathbb{R}^{dn}$ by
\begin{equation}\label{Rescale_B_Bar_s}
\overline{\mathbf{B}}(s) :=\frac{1}{\sqrt{\beta}}\mathbf{B}(\beta s)= \frac{1}{\sqrt{\beta}}\mathbf{W}(\beta s)-\frac{s}{\sqrt{\beta}}\mathbf{W}(\beta)\,, 
\end{equation}
where $\mathbf{W}:[0,\beta]\times\Omega\to\mathbb{R}^{dn}$ is the standard Wiener process in $\mathbb{R}^{dn}$. The process $\overline{\mathbf{B}}$ defined in \eqref{Rescale_B_Bar_s} is a standard Brownian bridge process, as verified\if\JOURNAL1
in Appendix~\ref{Appendix_section_b}. 
\fi
\if\JOURNAL2
in ~\ref{Appendix_section_b}. 
\fi
Using this re-scaling in the time variable, we have the alternative expression for the Brownian bridge path $\mathbf{x}_k^{\ell}(t)$,
\begin{equation}\label{scaling_BB_kell}
\mathbf x_{k}^{\ell}(t)=\sqrt{\beta}\, \overline{\mathbf{B}}_k(\frac{t}{\beta}) + (1-\frac{t}{\beta}) \mathbf x_k(0)+ \frac{t}{\beta} \mathbf x_{\ell}(0)\,,\ t\in[0,\beta]\,.
\end{equation}
To simplify the notation,  we define for all indices $k$ and ${\ell}$
\[ 
\begin{aligned}
a_{k\ell}(\beta,\mathbf{x}) &:= e^{-\frac{|\mathbf x_k(0)-\mathbf x_{\ell}(0)|^2}{2\beta}}\,,\\
\gamma_{k\ell}(\beta,\mathbf{x}) &:= \int_0^\beta \Tilde{V}_{k\ell}(\mathbf{x}(t))\,{\mathrm{d}} t\,,\\
\textup{and\ }\ b_{k\ell}(\beta,\mathbf{x}) &:= e^{-\gamma_{k\ell}(\beta,\mathbf{x}) }\,,\\
\end{aligned}
\]
then the matrix element $\mathcal{W}_{k\ell}(\beta,\mathbf{x})$ in \eqref{W_mat_element_no_nuclei} and its derivative with respect to $\beta$ can be written as
\[ 
    \mathcal{W}_{k\ell}(\beta,\mathbf{x}) = a_{k\ell}(\beta,\mathbf{x})\,b_{k\ell}(\beta,\mathbf{x})\,,
\]
and
\begin{equation}\label{derivative_W_kell}
\partial_\beta\, \mathcal{W}_{k\ell}(\beta,\mathbf{x})=\frac{\partial a_{k\ell}(\beta,\mathbf{x})}{\partial\beta}\,b_{k\ell}(\beta,\mathbf{x}) + a_{k\ell}(\beta,\mathbf{x})\,\frac{\partial b_{k\ell}(\beta,\mathbf{x})}{\partial\beta}\,.
\end{equation}
Specifically, we have
\begin{equation} \label{derivative_a_and_b}
\begin{aligned}
\frac{\partial a_{k\ell}(\beta,\mathbf{x})}{\partial\beta} &= \exp{\big(-\frac{|\mathbf x_k(0)-\mathbf x_{\ell}(0)|^2}{2\beta}\big)}\,\frac{|\mathbf x_k(0)-\mathbf x_{\ell}(0)|^2}{2\beta^2}\,,\\
 \frac{\partial b_{k\ell}(\beta,\mathbf{x})}{\partial\beta} &= e^{-\gamma_{k\ell}(\beta,\mathbf{x}) } \, \Big(-\frac{\partial }{\partial\beta}\,\gamma_{k\ell}(\beta,\mathbf{x})\Big)\,. \\
 \end{aligned}
 \end{equation} 
 Using the standard Brownian bridge process $\overline{\mathbf{B}}(s)$ with time variable $s\in[0,1]$ and introducing the shorthand notation
 \[  \Bar{\mathbf{x}}_{k}^{\ell}(s) := \sqrt{\beta}\, \overline{\mathbf{B}}_k(s) + (1-s)\, \mathbf x_k(0)+ s\,\mathbf x_{\ell}(0)\,,\textup{ for }s\in[0,1]\, ,\]
 we have the expression
 \[  
 \gamma_{k\ell}(\beta,\mathbf{x})=\int_0^1 \Big( \frac{1}{2}\big| \Bar{\mathbf{x}}_{k}^{\ell}(s) \big|^2 + \sum_{j=1,j\ne k}^n \frac{\lambda}{2}\, \frac{1}{\big| \Bar{\mathbf{x}}_{k}^{\ell}(s) - \Bar{\mathbf{x}}_{j}^{\nu_j}(s) \big|} \Big) \beta\,{\mathrm{d}} s=:\beta\int_0^1\tilde{V}_{k\ell}\big(\xbar(s)\big)\mathrm{d}s.
 \]
 Thus by employing the rescaling technique \eqref{scaling_BB_kell}, we avoid the derivative $\partial_\beta \mathbf{B}(t)$ and obtain
 \[  
\partial_\beta \gamma_{k\ell}(\beta,\mathbf{x})=\int_0^1\Big( \Tilde{V}_{k\ell}\big({\Bar{\mathbf{x}}(s)}\big)+\frac{\sqrt{\beta}}{2}\,\overline{\mathbf{B}}(s)\cdot\nabla\Tilde{V}_{k\ell}\big({\Bar{\mathbf{x}}(s)}\big) \Big)\,\mathrm{d}s.
 \]
 A more detailed expression for the implementation with pairwise Coulomb repulsion and a confining harmonic potential is provided\if\JOURNAL1
in Appendix~\ref{Appendix_section_a}.
 \fi
 \if\JOURNAL2
in ~\ref{Appendix_section_a}.
 \fi
 
\smallskip
With the absence of nuclei, the partition function $\mathcal{Z}_{\nu}(\beta)$ becomes
\[  
\mathcal{Z}_{\nu}(\beta) = \int_{\mathbb{R}^{3n}}\mathbb{E}\Big[ \frac{\mathrm{det}(\mathcal{W}(\mathbf{x})}{n!(2\pi\beta)^{3n/2}} \Big]\,\mathrm{d}\mathbf{x}\,,
\]
and the derivative of $\mathcal{Z}_{\nu}(\beta)$ with respect to $\beta$ is given by
\[ 
\begin{aligned}
\partial_\beta\, \mathcal{Z}_{\nu}(\beta) &= \int_{\mathbb{R}^{3n}}\frac{1}{n!}\mathbb{E}\Big[ \frac{1}{(2\pi\beta)^{3n/2}}\frac{\partial}{\partial\beta}\mathrm{det}\big(\mathcal{W}(\mathbf{x})\big)
-\frac{3n}{2\beta(2\pi)^{3n/2}}\beta^{-\frac{3n+2}{2}}%
\mathrm{det}\big( \mathcal{W}(\mathbf{x}) \big) \Big]\,\mathrm{d}\mathbf{x}\\
&=  \int_{\mathbb{R}^{3n}} \frac{1}{n!}\frac{1}{(2\pi\beta)^{3n/2}}\mathbb{E}\Big[
\mathrm{Tr}\Big( \mathrm{adj}\big(\mathcal{W}(\mathbf{x})\big)\, \frac{\partial \mathcal{W}(\mathbf{x}) }{\partial \beta}\Big)-\frac{3n}{2\beta}\,\mathrm{det}\big(\mathcal{W}(\mathbf{x})\big) \Big]\,\mathrm{d}\mathbf{x}\,.%
\end{aligned} 
\]
Therefore, based on the approximation of the partition function $\mathcal{Z}_{\nu}$, the system mean-field energy $h$ can be approximated by 
\begin{equation}\label{mean-field h approx. formula}
\begin{aligned}
h_{\nu} &:= -\frac{\partial_\beta\, \mathcal{Z}_{\nu}(\beta)}{\mathcal{Z}_{\nu}(\beta)}\\
&%
= -\frac{\int_{\mathbb{R}^{3n}} \mathbb{E}\Big[ \mathrm{Tr}\Big( \mathrm{adj}\big(\mathcal{W}(\mathbf{x})\big)\, \frac{\partial \mathcal{W}(\mathbf{x}) }{\partial \beta}\Big)-\frac{3n}{2\beta}\,\mathrm{det}\big(\mathcal{W}(\mathbf{x})\big) \Big]\,\mathrm{d}\mathbf{x_0}}{\int_{\mathbb{R}^{3n}} \mathbb{E}\Big[ \mathrm{det}\big(\mathcal{W}(\mathbf{x})\big) \Big]\,\mathrm{d}\mathbf{x_0}}\,,
\end{aligned}
\end{equation}
where $\frac{\partial \mathcal{W}(\mathbf{x}) }{\partial \beta}$ is computed by the element-wise formula \eqref{derivative_W_kell}.

In the numerical experiments, the integrals in the $\mathbb{R}^{3n}$ space appearing in \eqref{mean-field h approx. formula} are estimated with the Monte Carlo integration method by independently sampling initial states $\mathbf{x}(0\,;i)\in \mathbb{R}^{3n}$ for $i=1,2,\cdots,M_{x}$, and sampling, for each initial state $\mathbf{x}(0\,;i)$, an associated independent Brownian bridge path $\mathbf{B}(t\,;i)$ in $\mathbb{R}^{3n}$, which gives a sample for the whole path $\mathbf{x}(t\,;i)$ for $t\in[0,\beta]$ following \eqref{x_t_BB_path}. For simplicity of notation, we omit the variable $t$ for the stochastic Brownian paths $\mathbf{x}$ in \eqref{mean-field h approx. formula}.
More details on statistical analysis of simulations and  estimators of the mean-field approximation $\Bar{h}_v$ are provided\if\JOURNAL1
in Appendix~\ref{Appendix_section_c}.
\fi
\if\JOURNAL2
in ~\ref{Appendix_section_c}.
\fi

Similarly to Section~\ref{Subsection:Fermi_gas}, our sampling of Brownian bridge paths is based on a discretization scheme in time with step size $\Delta t$, and for evaluating the element-wise derivative $\frac{\partial \mathcal{W}(\mathbf{x}) }{\partial \beta}$,
we use the same numerical integration scheme with trapezoidal rule as shown in \eqref{trace_det2}, which gives a discretization error proportional to $\Delta t^2$ in our result. Our numerical approximation of mean-field energy involving also the time discretization error is denoted by $\Bar{h}_\nu$.

By varying the time-step size $\Delta t$ and taking a large sample size $M_{x}=2^{26}$ for Monte Carlo integral evaluation, we obtain the approximate mean-field energy $\Bar{h}_{\nu}$ for a system containing $n=6$ electrons with Coulomb interaction under a confining harmonic trap $V_{\textup{HO}}(\mathbf{x})=\frac{1}{2}|\mathbf{x}|^2$ at varying inverse temperature $\beta$, and make a comparison with the reference result $h_{\textup{ref}}$ obtained by the established algorithm in \cite{dornheim_1}. The corresponding results are summarized in 
Table~\ref{table_compare_Dornheim_n=6_dim=3}. As the system inverse temperature increases from $\beta=0.5$ to $\beta=2$, the relative difference level in our approximation $\Bar{h}_\nu$ rises from $10^{-4}$ to $10^{-2}$.
\begin{table}[h]
\begin{center}
\begin{tabular}{ c c c c c c }
\hline\hline\\ [-2.3ex]\multicolumn{6}{c}{\ }\\[-2.3ex]
 $\beta$ & $h_{\textup{ref}}$ &  $\Delta t$ & $\Bar{h}_{\nu}$ & Relative diff. & Relative CI  \\ \multicolumn{6}{c}{\ }\\[-2.3ex]
 \hline
\multicolumn{6}{c}{\ } \\ [-1.5ex]
\multirow{2}{*}{$0.5$} & \multirow{2}{*}{$41.66(1)$} & $0.025$ & $41.655(3)$ & $1.14 \,\mathrm{e}(-4)$ & $6.48 \,\mathrm{e}(-5)$\\
 &  & $0.0125$ & $41.653(3)$ & $1.64 \,\mathrm{e}(-4)$ & $6.47 \,\mathrm{e}(-5)$\\
 [-1.5ex]\multicolumn{6}{c}{\ } \\
\multirow{2}{*}{$1.0$} & \multirow{2}{*}{$26.692(9)$}  & $0.025$ & $26.711(4)$ & $7.12 \,\mathrm{e}(-4)$ & $1.48 \,\mathrm{e}(-4)$\\
 &  & $0.0125$ & $26.716(4)$ & $8.87 \,\mathrm{e}(-4)$ & $1.38 \,\mathrm{e}(-4)$\\
[-1.5ex]\multicolumn{6}{c}{\ } \\
\multirow{2}{*}{$1.5$} & \multirow{2}{*}{$22.63(7)$} & $0.025$ & $22.774(8)$ & $6.37 \,\mathrm{e}(-3)$ & $3.34 \,\mathrm{e}(-4)$\\
 &  & $0.0125$ & $22.774(7)$ & $6.37 \,\mathrm{e}(-3)$ & $3.23 \,\mathrm{e}(-4)$\\
[-1.5ex]\multicolumn{6}{c}{\ }\\
\multirow{2}{*}{$2.0$} & \multirow{2}{*}{$22.1(5)$} & $0.025$ & $20.66(5)$ & $6.52 \,\mathrm{e}(-2)$ & $2.47 \,\mathrm{e}(-3)$\\
 &  &  $0.0125$ & $20.69(6)$ & $6.39 \,\mathrm{e}(-2)$ & $2.69 \,\mathrm{e}(-3)$\\
[-1.75ex]
\multicolumn{6}{c}{\ } \\ \hline\hline
\end{tabular} 
\end{center}
\caption{  %
Case \ref{Case_V2}: The approximate mean-field energy $\Bar{h}_{\nu}$ obtained by Monte Carlo approximation of \eqref{mean-field h approx. formula} with sample size $M_x=2^{26}$ for $d=3$, $n=6$, compared to the reference value from \cite{dornheim_1}. The digit in the parenthesis of $\Bar{h}_\nu $ values denotes the statistical error from the Monte Carlo integration method. The relative difference between $\Bar{h}_\nu$ and $h_{\mathrm{ref}}$ is computed by $|\Bar{h}_\nu-h_{\mathrm{ref}}|/h_{\mathrm{ref}}$, and the last column records the 95\% relative confidence interval of $\Bar{h}_\nu$.
}
\label{table_compare_Dornheim_n=6_dim=3}
\end{table}

Moreover, we implement our path integral method with Feymann-Kac formulation on a $2$-dimensional quantum dot model, considering pair-wise Coulomb repulsion between electrons and still the confining harmonic trap $V_{\mathrm{HO}}$. The system inverse temperature is fixed to be $\beta=0.3$ and $\beta=1$, respectively, while the number of electrons increases from $n=6$ to $20$ and from $n=3$ to $10$. This $2$-dimensional model is also studied in \cite{dornheim_1} with the reference values for the mean-field energy $h_{\mathrm{ref}}$ for benchmarking our approximation $\Bar{h}_\nu$. We summarize our numerical results obtained with different time step sizes $\Delta t$ 
in Table~\ref{table_compare_Dornheim_dim=2_fix_beta_change_n}.
\begin{table}[h]
\begin{center}
\begin{tabular}{ c c c c c c c }
\hline\hline\\ [-2.3ex]\multicolumn{7}{c}{\ }\\[-2.3ex]
 $\beta$ & $n$ & $h_{\textup{ref}}$ &  $\Delta t$ & $\Bar{h}_{\nu}$ & Relative diff. & Relative CI  \\ \multicolumn{7}{c}{\ }\\[-2.3ex]
 \hline
\multicolumn{7}{c}{\ } \\ [-1.5ex]
\multirow{7}{*}{$1$} & \multirow{2}{*}{$3$} & \multirow{2}{*}{$8.719(3)$} & $0.025$ & $8.717(3)$ & $2.74 \,\mathrm{e}(-4)$ & $3.83\,\mathrm{e}(-4)$\\
 &  &  & $0.0125$ & $8.717(3)$ & $2.75 \,\mathrm{e}(-4)$ & $3.95 \,\mathrm{e}(-4)$\\
 [-1.5ex]\multicolumn{7}{c}{\ } \\
 & \multirow{2}{*}{$6$} & \multirow{2}{*}{$22.82(5)$}  & $0.025$ & $22.82(3)$ & $2.13 \,\mathrm{e}(-4)$ & $1.20 \,\mathrm{e}(-3)$\\
 &  &  & $0.0125$ & $22.79(2)$ & $1.03 \,\mathrm{e}(-3)$ & $1.04 \,\mathrm{e}(-3)$\\
[-1.5ex]\multicolumn{7}{c}{\ } \\
& \multirow{2}{*}{$10$} & \multirow{2}{*}{$49(3)$} & $0.025$ & $48.3(5)$ & $1.47\,\mathrm{e}(-2)$ & $9.92 \,\mathrm{e}(-3)$\\
&  &  & $0.0125$ & $48.5(4)$ & $1.04 \,\mathrm{e}(-2)$ & $8.89 \,\mathrm{e}(-3)$\\
 [-1.5ex]\multicolumn{7}{c}{\ }\\
\hline
\\[-4ex]\multicolumn{7}{c}{\ } \\
\multirow{7}{*}{$0.3$} & \multirow{2}{*}{$6$} & \multirow{2}{*}{$46.45(1)$} & $0.025$ & $46.44(1)$ & $2.17 \,\mathrm{e}(-4)$ & $2.83 \,\mathrm{e}(-4)$\\
 &  &  &  $0.0125$ & $46.45(1)$ & $9.79 \,\mathrm{e}(-6)$ & $2.85 \,\mathrm{e}(-4)$\\
 [-1.5ex]\multicolumn{7}{c}{\ }\\
 & \multirow{2}{*}{$10$} & \multirow{2}{*}{$84.92(4)$} & $0.025$ & $84.90(2)$ & $2.69\,\mathrm{e}(-4)$ & $2.85 \,\mathrm{e}(-4)$\\
 &  &  &  $0.0125$ & $84.89(2)$ & $3.04 \,\mathrm{e}(-4)$ & $2.86 \,\mathrm{e}(-4)$\\
 [-1.5ex]\multicolumn{7}{c}{\ }\\
  & \multirow{2}{*}{$20$} & \multirow{2}{*}{$203(1)$} & $0.025$ & $203.5(2)$ & $2.63\,\mathrm{e}(-3)$ & $1.04 \,\mathrm{e}(-3)$\\
 &  &  &  $0.0125$ & $203.4(2)$ & $2.07 \,\mathrm{e}(-3)$ & $8.85 \,\mathrm{e}(-4)$\\
[-1.75ex]
\multicolumn{7}{c}{\ } \\ \hline\hline
\end{tabular} 
\end{center}
\caption{  %
Case \ref{Case_V2}: The approximate mean-field energy $\Bar{h}_{\nu}$ obtained with equation \eqref{mean-field h approx. formula} using sample size $M_x=2^{22}$ for $d=2$, compared with the reference value from \cite{dornheim_1}. The digit in the parenthesis of $\Bar{h}_\nu$ data denotes the statistical error from the Monte Carlo integration method.
The relative difference between $\Bar{h}_\nu$ and $h_{\mathrm{ref}}$ is computed by $|\Bar{h}_\nu-h_{\mathrm{ref}}|/h_{\mathrm{ref}}$, and the last column records the 95\% relative confidence interval of $\Bar{h}_\nu$.
}
\label{table_compare_Dornheim_dim=2_fix_beta_change_n}
\end{table}
For both test cases 
with $\beta=1$ and $\beta=0.3$, our approximation $\Bar{h}_{\nu}$ employs $M_x=2^{22}$ samples for Monte Carlo integral evaluation and agrees with the reference value $h_{\mathrm{ref}}$ on their corresponding level of statistical confidence intervals. %

To further survey the bias level of our mean-field energy approximation $\Bar{h}_{\nu}$ in the test case \ref{Case_V2}, we compare our results with both the reference value $h_{\mathrm{ref}}$ in \cite{dornheim_1}, and with a Monte Carlo approximation $h_{\mathrm{tensor}}$ of the exact mean-field energy based on the tensor formulation \eqref{tensor}, for $n=6$ fermions in dimension $d=3$. The sample size $M_x$ employed for the Monte Carlo integral evaluation of $\Bar{h}_\nu$ is $2^{26}$. For the Monte Carlo integral approximation of the more computationally demanding tensor formula, however, the sample size is taken to be $2^{22}$, which leads to relatively large confidence intervals for the obtained $h_{\mathrm{tensor}}$ data. The numerical results for $h_{\mathrm{ref}}$, $h_{\mathrm{tensor}}$, and $\Bar{h}_{\nu}$ at different $\beta$ values are summarized in Table \ref{table_compare_htensor_and_Dornheim_n=6_dim=3}, from which we observe that the relative difference $|\Bar{h}_\nu-h_{\mathrm{tensor}}|/h_{\mathrm{tensor}}$ increases as $\beta$ increases. Also, the differences $|\Bar{h}_\nu-h_{\mathrm{tensor}}|$ and $|h_{\mathrm{ref}}-h_{\mathrm{tensor}}|$ are smaller than the statistical uncertainty of $h_{\mathrm{tensor}}$, suggesting that the statistical error from the high-dimensional Monte Carlo integral evaluation is still dominating.

\begin{table}[h]
\begin{center}
\begin{tabular}{ c c c c c c c }
\hline\hline\\ [-2.3ex]\multicolumn{6}{c}{\ }\\[-2.3ex]
 $\beta$ & $h_{\textup{ref}}$ &  $\Delta t$ & $h_{\mathrm{tensor}}$ & $\Bar{h}_{\nu}$ & %
 $\frac{|h_{\mathrm{tensor}}-\Bar{h}_{\nu}|}{h_{\mathrm{tensor}}}$
 & Relative CI of $h_{\mathrm{tensor}}$  \\ \multicolumn{7}{c}{\ }\\[-2.3ex]
 \hline
\multicolumn{7}{c}{\ } \\ [-1.5ex]
\multirow{2}{*}{$0.5$} & \multirow{2}{*}{$41.66(1)$} & $0.025$ & $41.65(7)$ & $41.655(3)$ & $7.20\, \mathrm{e}(-5)$ & $1.64 \,\mathrm{e}(-3)$\\
 & & $0.0125$ & $41.66(7)$ & $41.653(3)$ & $3.87 \,\mathrm{e}(-5)$ & $1.64 \,\mathrm{e}(-3)$\\
 [-1.5ex]\multicolumn{6}{c}{\ } \\
\multirow{2}{*}{$1.0$} & \multirow{2}{*}{$26.692(9)$}  & $0.025$ & $26.7(1)$ & $26.711(4)$ & $4.12 \,\mathrm{e}(-4)$ & $5.71 \,\mathrm{e}(-3)$\\
 &  & $0.0125$ & $26.7(1)$ & $26.716(4)$ & $5.99 \,\mathrm{e}(-4)$ & $5.24 \,\mathrm{e}(-3)$\\
[-1.5ex]\multicolumn{6}{c}{\ } \\
\multirow{2}{*}{$1.5$} & \multirow{2}{*}{$22.63(7)$} & $0.025$ & $22.8(2)$ & $22.774(8)$ & $1.13 \,\mathrm{e}(-3)$ & $1.05\,\mathrm{e}(-2)$\\
 &  & $0.0125$ & $22.8(2)$ & $22.774(7)$ & $1.09 \,\mathrm{e}(-3)$ & $8.59 \,\mathrm{e}(-3)$\\
[-1.75ex]
\multicolumn{6}{c}{\ } \\ \hline\hline
\end{tabular} 
\end{center}
\caption{  %
Case \ref{Case_V2}: The approximate mean-field energy $\Bar{h}_{\nu}$ obtained by Monte Carlo approximation of \eqref{mean-field h approx. formula} with sample size $M_x=2^{26}$ for $d=3$, $n=6$, compared to the reference value from \cite{dornheim_1} and from the tensor formula \eqref{tensor} with sample size $2^{22}$. The digit in the parenthesis of $h_{\mathrm{tensor}}$ and $\Bar{h}_\nu $ values denotes the statistical error from the Monte Carlo integration method. The relative difference between $\Bar{h}_\nu$ and $h_{\mathrm{tensor}}$ is computed by $|\Bar{h}_\nu-h_{\mathrm{tensor}}|/h_{\mathrm{tensor}}$, and the last column records the 95\% relative confidence interval of $h_{\mathrm{tensor}}$.
}
\label{table_compare_htensor_and_Dornheim_n=6_dim=3}
\end{table}

In addition, we implement the perturbation of the determinant formulation following \eqref{xi_approx}, \eqref{Z-perturbed} and \eqref{W_bar}, with the aim of obtaining a sensitivity study for the $\Bar{h}_\nu$ approximation. We compute first the approximated partition function $\Bar{Z}_{\mathrm{perturb}}(\beta)$, using the expected values of $\mathrm{det}\big(\widetilde{\mathcal{W}}\big)$ based on 100 independent samples of perturbed Brownian paths, and then apply a central difference quotient formula to approximate the mean-field value $\Bar{h}_{\mathrm{perturb}}(\beta)$ as %
\begin{equation}\label{h-perturbed} 
\Bar{h}_{\mathrm{perturb}}(\beta):= \frac{\Bar{Z}_{\mathrm{perturb}}(\beta+\Delta\beta)-\Bar{Z}_{\mathrm{perturb}}(\beta-\Delta\beta)}{2\Delta\beta}.%
\end{equation}
The results of $\Bar{h}_{\mathrm{perturb}}$ and the non-perturbed approximation $\Bar{h}_\nu$ using the determinant formulation are summarized in Table \ref{Table_compare_perturb_BB_with_W_matrix}. These results are then compared with the reference mean-field value $h_{\mathrm{ref}}$. Specifically, for $n=3$, we obtain $h_{\mathrm{ref}}$ using the tensor formula \eqref{tensor}, while for $n=6$, we rely on the result from \cite{dornheim_1} to obtain a reliable reference value $h_{\mathrm{ref}}$. 
The relative error of the perturbed mean-field value $\Bar{h}_{\mathrm{perturb}}$ is at most one order of magnitude higher than the relative error of $\Bar{h}_\nu$ derived from the determinant formulation \eqref{mean-field h approx. formula}, which indicates that $\Bar{h}_{\mathrm{perturb}}$ can be used to numerically roughly estimate the accuracy of $\bar h_\nu$ without using a reference value. 
In both test cases with $\beta=1.0$ and $\beta=1.5$, we employ $M_x=2^{22}$ samples for the Monte Carlo evaluation of the integral in the high-dimensional $\mathbb{R}^{3n}$ space, and for the central difference formula \eqref{h-perturbed} we use $\Delta \beta = 0.01$.
\begin{table}[h]
\begin{center}
\begin{tabular}{ c c c c c c c c }
\hline\hline\\ [-2.3ex]\multicolumn{8}{c}{\ }\\[-2.3ex]
 $\beta$ & $n$ & $h_{\textup{ref}}$ &  $\Delta t$ & $\Bar{h}_{\nu}$ & $\frac{|h_{\mathrm{ref}}-\Bar{h}_{\nu}|}{h_{\mathrm{ref}}}$ & $\Bar{h}_{\mathrm{perturb}}$ & $\frac{|h_{\mathrm{ref}}-\Bar{h}_{\mathrm{perturb}}|}{h_{\mathrm{ref}}}$  \\ \multicolumn{8}{c}{\ }\\[-2.3ex]
 \hline
\multicolumn{8}{c}{\ } \\ [-1.5ex]
\multirow{4}{*}{$1.0$} & \multirow{2}{*}{$3$} & \multirow{2}{*}{$11.355(3)$} &  $0.025$ & $11.356(3)$  & $1.25\,\mathrm{e}(-4)$ & $11.337(4)$ & $1.55\,\mathrm{e}(-3)$ \\
 &  &  & $0.0125$ & $11.356(3)$ & $7.65\,\mathrm{e}(-5)$ & $11.337(3)$ & $1.63\,\mathrm{e}(-3)$ \\
 [-1.5ex]\multicolumn{8}{c}{\ } \\
 & \multirow{2}{*}{$6$} & \multirow{2}{*}{$26.692(9)$}&  $0.025$ & $26.72(1)$ & $1.52\,\mathrm{e}(-4)$ & $26.60(1)$ & $3.40\,\mathrm{e}(-3)$\\
 &  &  & $0.0125$ & $26.71(1)$ & $1.08\,\mathrm{e}(-4)$ & $26.60(1)$ & $3.56\,\mathrm{e}(-3)$\\
 [-1.5ex]\multicolumn{8}{c}{\ }\\
\hline
\\[-4ex]\multicolumn{8}{c}{\ } \\

\multirow{4}{*}{$1.5$} & \multirow{2}{*}{$3$} & \multirow{2}{*}{$9.157(2)$} &  $0.025$ & $9.163(4)$ & $6.23 \,\mathrm{e}(-4)$ & $9.129(5)$ & $3.02 \,\mathrm{e}(-3)$\\
 &  &  &  $0.0125$ & $9.159(4)$ & $1.86 \,\mathrm{e}(-4)$ & $9.132(5)$ & $2.71 \,\mathrm{e}(-3)$\\
 [-1.5ex]\multicolumn{8}{c}{\ }\\
 & \multirow{2}{*}{$6$} & \multirow{2}{*}{$22.63(7)$} &  $0.025$ & $22.83(3)$ & $8.79\,\mathrm{e}(-3)$ & $22.60(5)$ & $1.20 \,\mathrm{e}(-3)$\\
 &  &  &  $0.0125$ & $22.85(2)$ & $9.87 \,\mathrm{e}(-3)$ & $22.55(5)$ & $3.55 \,\mathrm{e}(-3)$\\
[-1.75ex]
\multicolumn{8}{c}{\ }
 \\ \hline\hline
\end{tabular} 
\end{center}
\caption{  %
Case \ref{Case_V2}: The approximate mean-field energy $\Bar{h}_{\nu}$ for $d=3$, obtained with equation \eqref{mean-field h approx. formula} and the corresponding result $\Bar{h}_{\mathrm{perturb}}$ obtained with the perturbation determinant formula \eqref{h-perturbed} using sample size $M_x=2^{22}$, compared with the reference value $h_{\mathrm{ref}}$ obtained by the tensor formula based on \eqref{tensor} and from the benchmark result in \cite{dornheim_1}. The digit in the parenthesis of $\Bar{h}_\nu$ and $\Bar{h}_{\mathrm{perturb}}$ denotes the statistical error from the Monte Carlo integration method. All the values in this table are computed with parameter $c_*=2$, which is determined according to the cycle length for the given number of particles $n$. Numerical tests with $c_\ast$ varying from 0.6 to 2 indicate an insensitivity with respect to this parameter, since $|\Bar{h}_{\mathrm{perturb}}-h_{\mathrm{ref}}|/h_{\mathrm{ref}}$ remains smaller than $5\,\mathrm{e}(-3)$, which is within the range of the statistical error.
}
\label{Table_compare_perturb_BB_with_W_matrix}
\end{table}

\section{Conclusions} \label{sec_conclusion}

We develop and analyze a new computational method for approximating the mean-field Hamiltonian of molecular dynamics suitable for estimation of canonical quantum observables in nuclei-electrons systems.
The method  is based on approximation of path integrals for fermions 
using Monte Carlo estimation of particular determinants on Brownian bridge paths.  In the presented numerical benchmarks we observe that
the derived determinant formulation is surprisingly accurate.
More specifically, the reported test examples show an approximation 
error at 
the same level as the statistical error in the reference solutions taken from \cite{dornheim_1}.
The approximation of mean-feld value $\Bar{h}_\nu$, based on discretizations of \eqref{mean-field h approx. formula}, is also consistent with the numerical results of $\Bar{h}_{\mathrm{perturb}}$, obtained from the sensitivity analysis \eqref{xi_approx}, \eqref{Z-perturbed}, \eqref{W_bar}, \eqref{h-perturbed}, shown in Table~\ref{Table_compare_perturb_BB_with_W_matrix}.
The perturbed mean-field value $\Bar{h}_{\mathrm{perturb}}$ thereby provides a computational error indicator for
$\bar h_\nu$ when no reference solution is available. 
It remains to  explain theoretically the accuracy of the approximation based on the formulation \eqref{Z_det} from the first principles. 
We remark that similar to other path dependent methods based on the Feynman-Kac formulation  the presented method exhibits an increasing bias as the  temperature decreases.

For the test case \ref{Case_V1} with non-interacting Fermi gas under an external harmonic oscillator potential $V_{\mathrm{HO}}(\mathbf{x})=\frac{1}{2}|\mathbf{x}|^2$, we compare  our numerical approximation of the partition function $\bar{Z}_n(\beta)$ and the mean-field energy $\Bar{h}(\beta)$ with their corresponding exact values. Since the potential $V_{\mathrm{HO}}$ is particle-wise separable, our formulation \eqref{trace_det2} is asymptotically exact as the Monte Carlo integration sample size $M_x\to \infty$ and the time step size $\Delta t\to 0$.  
From Tables~\ref{Table_approx_Z_HO}, \ref{Table_approx_mean-field_HO}, and Figure~\ref{Fig:rel_error_pure_HO}, we observe no biases in the relative error of the approximations as $M_x$ increases, and the statistical error 
from the Monte Carlo approximation of the high-dimensional integral in the Feynman-Kac  formulation \eqref{trace_det2} is dominating. Meanwhile, the time discretization error related to $\Delta t$ is smaller than the statistical uncertainty.

For the test case \ref{Case_V2}, which involves the pairwise Coulomb interactions between electrons and a confining harmonic trap, our numerical approximation to the mean-field energy $\Bar{h}_\nu$ achieves comparable accuracy to the reference result in \cite{dornheim_1} for both the two models with dimensions $d=3$ and $d=2$, as shown in the Tables~\ref{table_compare_Dornheim_n=6_dim=3} and ~\ref{table_compare_Dornheim_dim=2_fix_beta_change_n}. Specifically, as the inverse temperature $\beta$ and the number of fermions $n$ increases, the level of the statistical uncertainty for both methods rises. 
From Table~\ref{table_compare_Dornheim_n=6_dim=3} it is observed that the relative difference to the reference mean-field $|\Bar{h}_\nu - h_{\mathrm{ref}}|/h_{\mathrm{ref}}$ increases as $\beta$ increases, indicating an elevated bias level in the formulation \eqref{Z_det} for a lower system temperature.

In comparison with the reference result $h_{\mathrm{tensor}}$ using the permutation formula \eqref{tensor} we still observe 
from Table~\ref{table_compare_htensor_and_Dornheim_n=6_dim=3} that the statistical error from the Monte Carlo integration method is dominating, while for a system comprising $n$ electrons, the computational workload of the tensor formula $h_{\mathrm{tensor}}$ is $\mathcal{O}(n!)$, however, our approximation $\Bar{h}_\nu$ only requires a computational complexity proportional to $n^3$.

\section*{Acknowledgment}
This research was supported by
Swedish Research Council grant 2019-03725 and KAUST grant OSR-2019-CRG8-4033.3. 
The work of P.P. was supported in part by the U.S. Army Research Office Grant W911NF-19-1-0243. The computations were enabled by resources provided by the National Academic Infrastructure for Supercomputing in Sweden (NAISS) at PDC Center for High Performance Computing, KTH Royal Institute of Technology, partially funded by the Swedish Research Council through grant agreement no. 2022-06725.

\if\JOURNAL1
\bibliographystyle{abbrv}
\fi
\if\JOURNAL2
\bibliographystyle{elsarticle-num}
\fi

\appendix %
\section{Supplementary for Numerical Implementation}
\subsection{The derivative of the matrix element \texorpdfstring{$\mathcal{W}_{k\ell}(\beta,\mathbf{x})$}{Lg} with respect to \texorpdfstring{$\beta$}{Lg}}\label{Appendix_section_a}
\ \\  For a system comprising indistinguishable fermions with Coulomb repulsive interactions and under a confining harmonic oscillator potential, we have the formula for the matrix element $\mathcal{W}_{k\ell}(\beta,\mathbf{x})$, based on \eqref{W_mat_element_no_nuclei} and \eqref{V_tilde_nu_j}:
\[
\mathcal W_{k\ell}(\beta,\mathbf x) =
 e^{-\frac{|\mathbf x_k(0)-\mathbf x_{\ell}(0)|^2}{2\beta}} 
e^{ -\int_0^\beta \Tilde{V}_{k,\ell}(\mathbf{x}(t)){\mathrm{d}} t}\,,       
\]
where
\[
\Tilde{V}_{k,\ell}\big(\mathbf{x}(t)\big)=\frac{1}{2} |\mathbf{x}_k^{\ell}(t)|^2 + \sum_{j=1,j\ne k}^n \frac{\lambda}{2|\mathbf{x}_k^{\ell}(t) - \mathbf{x}_j^{\nu_j}(t)|}, 
\]
Using the scaling formula \eqref{scaling_BB_kell} we can first rewrite the expression for $\gamma_{k,\ell}(\beta,\mathbf{x})$ with a change of the integration variable as
\begin{equation}\label{gamma_kell_integral_with_s}
    \begin{split}
    \gamma_{k,\ell}(\beta,\mathbf{x}) &= \int_0^\beta \Tilde{V}_{k,\ell}(\mathbf{x}(t))\,{\mathrm{d}} t\\
    &=\int_0^\beta \frac{1}{2} |\mathbf{x}_k^{\ell}(t)|^2 + \sum_{j=1,j\ne k}^n \frac{\lambda}{2|\mathbf{x}_k^{\ell}(t) - \mathbf{x}_j^{\nu_j}(t)|} \,\mathrm{d} t\\
    & \dots \big\{ \ \textup{Plugging in the scaling formulae \eqref{scaling_BB_kell}\,} \ \big\}\\
    &= \int_0^\beta  \frac{1}{2}\, \big|\underbrace{\sqrt{\beta}\, \Bar{\mathbf{B}}_k(\frac{t}{\beta}) + (1-\frac{t}{\beta}) \mathbf x_k(0)+ \frac{t}{\beta} \mathbf x_{\ell}(0)}_{=:\Bar{\mathbf{x}}_{k}^{\ell}(\frac{t}{\beta}) }\big|^2 \,\mathrm{d}t  \\
  \qquad\;\;\; &+ \int_0^\beta \sum_{j=1,j\ne k}^n \frac{\lambda}{2}\frac{\mathrm{d}t}{\big|\Bar{\mathbf{x}}_{k}^{\ell}(\frac{t}{\beta}) - \Bar{\mathbf{x}}_{j}^{\nu_j}(\frac{t}{\beta}) \big|}  \,\\
    & \dots \big\{ \ \textup{Taking the change of variable } \,s:=\frac{t}{\beta},\ \mathrm{d}t = \beta\,\mathrm{d}s \ \big\}\\
        &=\quad \int_0^1  \frac{1}{2}\, \big|\Bar{\mathbf{x}}_{k}^{\ell}(s)\big|^2\beta \,\mathrm{d}s  
         + \int_0^1 \sum_{j=1,j\ne k}^n \frac{\lambda}{2}\frac{1}{\big|\Bar{\mathbf{x}}_{k}^{\ell}(s) - \Bar{\mathbf{x}}_{j}^{\nu_j}(s) \big|}  \,\beta\,\mathrm{d}s.\\
    \end{split}
\end{equation}
We have 
\begin{equation}\label{derivative_x_bar_kell}
\frac{\partial}{\partial \beta}\big| \Bar{\mathbf{x}}_{k}^{\ell}(s) \big|^2 = 2\, \Bar{\mathbf{x}}_{k}^{\ell}(s) \cdot \frac{\partial}{\partial \beta} \Bar{\mathbf{x}}_{k}^{\ell}(s) = \Bar{\mathbf{x}}_{k}^{\ell}(s) \cdot \frac{1}{\sqrt{\beta}}\Bar{\mathbf{B}}_k(s)\, ,
\end{equation}
and
\begin{equation}\label{derivative_x_bar_j_nu_j}
\begin{split}
&\frac{\partial}{\partial \beta}\frac{1}{\big| \Bar{\mathbf{x}}_{k}^{\ell}(s) - \Bar{\mathbf{x}}_{j}^{\nu_j}(s) \big|} \\
 &= -\frac{1}{\big| \Bar{\mathbf{x}}_{k}^{\ell}(s) - \Bar{\mathbf{x}}_{j}^{\nu_j}(s) \big|^3}\, \big( \Bar{\mathbf{x}}_{k}^{\ell}(s) - \Bar{\mathbf{x}}_{j}^{\nu_j}(s) \big) \cdot \frac{\partial}{\partial \beta} \big( \Bar{\mathbf{x}}_{k}^{\ell}(s) - \Bar{\mathbf{x}}_{j}^{\nu_j}(s) \big)\\
&= -\frac{1}{\big| \Bar{\mathbf{x}}_{k}^{\ell}(s) - \Bar{\mathbf{x}}_{j}^{\nu_j}(s) \big|^3}\, \big( \Bar{\mathbf{x}}_{k}^{\ell}(s) - \Bar{\mathbf{x}}_{j}^{\nu_j}(s) \big) \cdot \Big( \frac{1}{2\sqrt{\beta}}\Bar{\mathbf{B}}_k(s) - \frac{1}{2\sqrt{\beta}}\Bar{\mathbf{B}}_j(s) \Big)\,.
\end{split}
\end{equation}
Based on equations \eqref{gamma_kell_integral_with_s}, \eqref{derivative_x_bar_kell}, and \eqref{derivative_x_bar_j_nu_j}, we further obtain
\begin{equation}\label{derivative_gamma}
\begin{aligned}
  \quad &  \frac{\partial }{\partial\beta}\,\gamma_{k,\ell}(\beta,\mathbf{x}) = \frac{\partial}{\partial\beta}\Big( \int_0^1 \big( \frac{1}{2}\big| \Bar{\mathbf{x}}_k^{\ell}(s) \big|^2 + \sum_{j=1,j\ne k}^n \frac{\lambda}{2}\, \frac{1}{\big| \Bar{\mathbf{x}}_{k}^{\ell}(s) - \Bar{\mathbf{x}}_{j}^{\nu_j}(s) \big|} \big) \beta\,{\mathrm{d}} s \Big)\\
 & \begin{aligned}
 =& \int_0^1 \big( \frac{1}{2}\,\frac{\partial}{\partial\beta}\big| \Bar{\mathbf{x}}_k^{\ell}(s) \big|^2 +\frac{\lambda}{2} \sum_{j=1,j\ne k}^n \frac{\partial}{\partial\beta}\frac{1}{\big| \Bar{\mathbf{x}}_{k}^{\ell}(s) - \Bar{\mathbf{x}}_{j}^{\nu_j}(s) \big|} \big) \beta\,{\mathrm{d}} s\\
 &+ \int_0^1 \big( \frac{1}{2}\big| \Bar{\mathbf{x}}_k^{\ell}(s) \big|^2 + \sum_{j=1,j\ne k}^n \frac{\lambda}{2}\, \frac{1}{\big| \Bar{\mathbf{x}}_{k}^{\ell}(s) - \Bar{\mathbf{x}}_{j}^{\nu_j}(s) \big|} \big) \,{\mathrm{d}} s\\
 \end{aligned}\\
 &  \begin{aligned}
 =& \int_0^1 \Big( \frac{1}{2}\,\Bar{\mathbf{x}}_{k}^{\ell}(s) \cdot \frac{1}{\sqrt{\beta}}\Bar{\mathbf{B}}_k(s) \\
  &-\frac{\lambda}{2} \sum_{j=1,j\ne k}^n \frac{\big( \Bar{\mathbf{x}}_{k}^{\ell}(s) - \Bar{\mathbf{x}}_{j}^{\nu_j}(s) \big)}{\big| \Bar{\mathbf{x}}_{k}^{\ell}(s) - \Bar{\mathbf{x}}_{j}^{\nu_j}(s) \big|^3} \cdot \frac{1}{2\sqrt{\beta}}\big(\Bar{\mathbf{B}}_k(s) - \Bar{\mathbf{B}}_j(s) \big) \Big) \beta\,{\mathrm{d}} s\\
 &+ \int_0^1 \Big( \frac{1}{2}\big| \Bar{\mathbf{x}}_k^{\ell}(s) \big|^2 + \sum_{j=1,j\ne k}^n \frac{\lambda}{2}\, \frac{1}{\big| \Bar{\mathbf{x}}_{k}^{\ell}(s) - \Bar{\mathbf{x}}_{j}^{\nu_j}(s) \big|} \Big) \,{\mathrm{d}} s\\
 \end{aligned}\\
 & \begin{aligned}
 =& \int_0^\beta \Big( \frac{1}{2}\,{\mathbf{x}}_{k}^{\ell}(t) \cdot \frac{1}{\beta}{\mathbf{B}}_k(t)\\
 &-\frac{\lambda}{2} \sum_{j=1,j\ne k}^n \frac{\big( {\mathbf{x}}_{k}^{\ell}(t) - {\mathbf{x}}_{j}^{\nu_j}(t) \big)}{\big| {\mathbf{x}}_{k}^{\ell}(t) - {\mathbf{x}}_{j}^{\nu_j}(t) \big|^3} \cdot \frac{1}{2\beta}\big({\mathbf{B}}_k(t) - {\mathbf{B}}_j(t) \big) \Big) \,{\mathrm{d}} t\\
 &+ \int_0^\beta \Big( \frac{1}{2}\big| {\mathbf{x}}_k^{\ell}(t) \big|^2 + \sum_{j=1,j\ne k}^n \frac{\lambda}{2}\, \frac{1}{\big| {\mathbf{x}}_{k}^{\ell}(t) - {\mathbf{x}}_{j}^{\nu_j}(t) \big|} \Big)\frac{1}{\beta} \,{\mathrm{d}} t\,.\\
 \end{aligned}\\
\end{aligned}
\end{equation}
Now combining equations \eqref{derivative_gamma}, \eqref{derivative_W_kell}, and \eqref{derivative_a_and_b}, we obtain at the derivative of the matrix element $\partial_\beta\,\mathcal{W}_{k,\ell}(\beta,\mathbf{x})$.
\subsection{The time-rescaling technique on the Brownian bridge process \texorpdfstring{$\mathbf{B}(t)$}{Lg} }\label{Appendix_section_b}
\ \newline
Applying the covariance function of the standard Wiener process $\mathbf{W}:\mathbb{R}^+\times\Omega\to\mathbb{R}$,
\[ \mathrm{Cov}\big(\mathbf{W}(\tau_1),\mathbf{W}(\tau_2)\big)=\min{(\tau_1,\tau_2)}\,, \]
we obtain the covariance function of $\overline{\mathbf{B}}(s)$ as defined in \eqref{Rescale_B_Bar_s} for two time values $\tau_1,\tau_2\in[0,1]$ as
\[
\begin{aligned}
    &\mathrm{Cov}\big(\overline{\mathbf{B}}(\tau_1),\overline{\mathbf{B}}(\tau_2)\big) \\
    &= \mathrm{Cov}\big( \frac{1}{\sqrt{\beta}}\mathbf{W}(\tau_1\beta)-\frac{\tau_1}{\sqrt{\beta}}\mathbf{W}(\beta) , \frac{1}{\sqrt{\beta}}\mathbf{W}(\tau_2\beta)-\frac{\tau_2}{\sqrt{\beta}}\mathbf{W}(\beta) \big)\\
    & = \frac{1}{\beta}\Big( \mathrm{Cov}\big( \mathbf{W}(\tau_1\beta) , \mathbf{W}(\tau_2\beta) \big) - \tau_1\mathrm{Cov}\big( \mathbf{W}(\tau_2\beta) , \mathbf{W}(\beta) \big) - \tau_2\mathrm{Cov}\big( \mathbf{W}(\tau_1\beta) , \mathbf{W}(\beta) \big)\Big) \\
    &\quad +\frac{1}{\beta}\, \tau_1\tau_2\mathrm{Cov}\big( \mathbf{W}(\beta)  , \mathbf{W}(\beta) \big)\\
    &= \frac{1}{\beta}\Big( \min{(\tau_1\beta,\tau_2\beta)}-\tau_1\tau_2\beta-\tau_2\tau_1\beta+\tau_1\tau_2\beta \Big)\\
    &= \min{(\tau_1,\tau_2)}-\tau_1\tau_2\,,
\end{aligned}
\]
which satisfies the definition of a standard Brownian bridge process defined with the time range $s\in[0,1]$.
The Brownian bridge $\mathbf{B}_k(t)$ in the equation \eqref{x_k_ell(t)_path_def_1} can thus be expressed with a scaling factor $\beta$ as
\begin{equation} \label{scaling_BB_k}
\mathbf{B}_k(t) \overset{\mathrm{d}}{=} \sqrt{\beta}\, \overline{\mathbf{B}}_k(\frac{t}{\beta})\,,\ t\in[0,\beta]\,, 
\end{equation}
where $\overline{\mathbf{B}}_k(s)$ denotes a standard Brownian bridge starting and ending at the same point $\mathbf{x}_k(0)$ with the time variable $s\in[0,1]$, and `$\mathrm{d}$' denotes equality of the two random processes in distribution. 
\section{ Statistical confidence interval of the approximation of mean-field \texorpdfstring{$\Bar{h}_\nu$}{Lg}} \label{Appendix_section_c}
Based on the approximation \eqref{Z_det} of the partition function $\mathcal{Z}_\nu$, we have for the mean-field energy the approximate formula \eqref{mean-field h approx. formula}
\[
h_{\nu} = -\frac{\partial_\beta\, \mathcal{Z}_{\nu}(\beta)}{\mathcal{Z}_{\nu}(\beta)}
= -\frac{\int_{\mathbb{R}^{3n}} \mathbb{E}\Big[ \mathrm{Tr}\Big( \mathrm{adj}\big(\mathcal{W}(\mathbf{x})\big)\, \frac{\partial \mathcal{W}(\mathbf{x}) }{\partial \beta}\Big)-\frac{3n}{2\beta}\,\mathrm{det}\big(\mathcal{W}(\mathbf{x})\big) \Big]\,\mathrm{d}\mathbf{x_0}}{\int_{\mathbb{R}^{3n}} \mathbb{E}\Big[ \mathrm{det}\big(\mathcal{W}(\mathbf{x})\big) \Big]\,\mathrm{d}\mathbf{x_0}}
=:\frac{\mathbb{E}[A]}{\mathbb{E}[B]}\,,
\]
where $A$ and $B$ are the random variables corresponding to the Monte Carlo approximation of the integration in the numerator and the denominator, respectively. For evaluating the Monte Carlo integrals, a sample set of size $N$ is utilized, which yields the approximation of the mean-filed $\Bar{h}_{\nu,N}$
\begin{equation}\label{h_nu_N}
    \Bar{h}_{\nu,N}=\frac{\frac{1}{N}\sum_{n=1}^N A_n}{\frac{1}{N}\sum_{n=1}^N B_n}\,,
\end{equation}
where $\{A_n\}_{n=1}^N$ and $\{B_n\}_{n=1}^N$ are sample sets based on the identically and independently drawn samples of paths $\{\mathbf{x}^{(n)}\}_{n=1}^N$ with initial position $\mathbf{x}^{(n)}_0$ following the mixed normal distribution with probability density $p(x)$, which is described by \eqref{p_x_density}. Specifically we have
\[
\begin{aligned}
    A_n&=\Big(\mathrm{Tr}\big( \mathrm{adj}\big(\mathcal{W}(\mathbf{x}^{(n)})\big)\, \frac{\partial \mathcal{W}(\mathbf{x}^{(n)}) }{\partial \beta}\big)-\frac{3n}{2\beta}\,\mathrm{det}\big(\mathcal{W}(\mathbf{x}^{(n)})\big)\Big)/p(\mathbf{x}^{(n)}_0)\,,\\
    B_n&=\mathrm{det}\big(\mathcal{W}(\mathbf{x}^{(n)})\big)/p(\mathbf{x}^{(n)}_0)\,.\\
\end{aligned}
\]
In this section, we derive estimates of the above statistical approximation of $h_\nu$ based on the Berry-Esseen theorem.

\begin{figure}[htbp]
    \centering

    \subfigure[$g_\epsilon(z)$]{\includegraphics[width=0.49\textwidth]{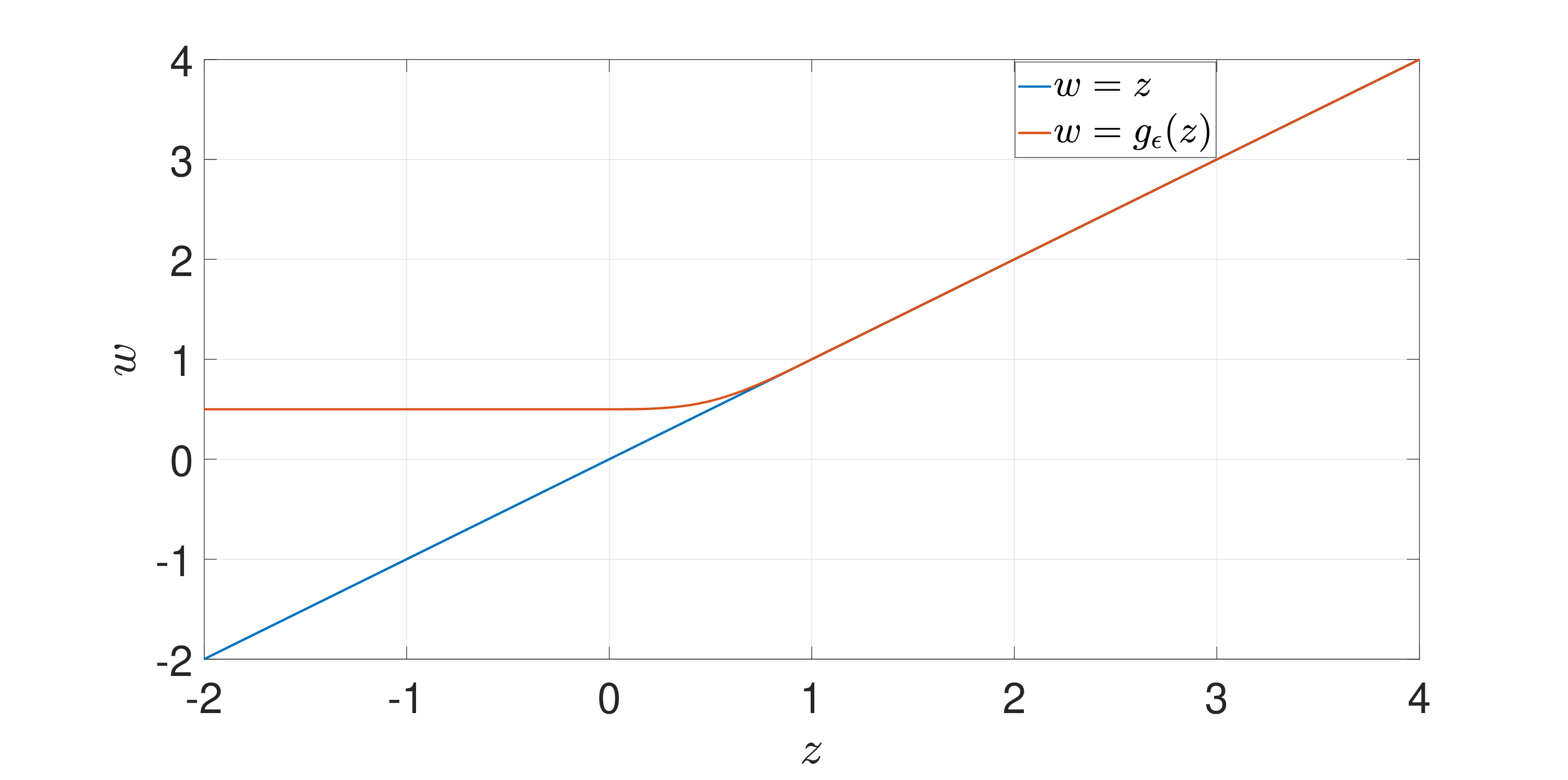}}
    \hfill
    \subfigure[$1/g_\epsilon(b+y)$]{\includegraphics[width=0.49\textwidth]{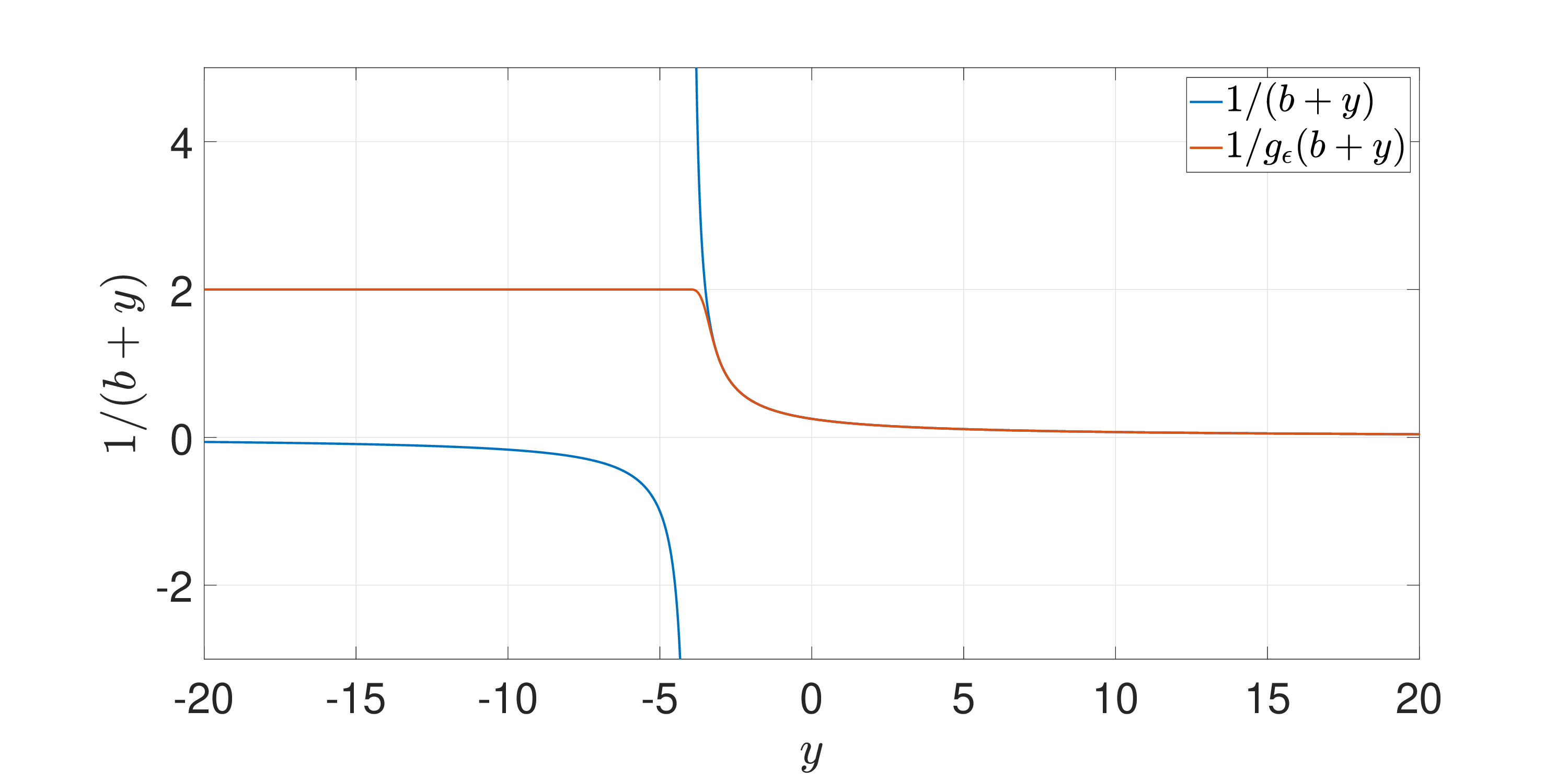}}

    \caption{Plot of the function $g_\epsilon(z)$ and a comparison of the two functions $1/(b+y)$ and $1/g_\epsilon(b+y)$, with parameters $\epsilon=0.5$, $b=4$. }
    \label{fig:plot_of_function_g}
\end{figure}

Applying the central limit theorem, letting
\[
\begin{aligned}
    \frac{1}{N}\sum_{n=1}^N A_n&=:a+\frac{\xi}{\sqrt{N}}\,,\\
    \frac{1}{N}\sum_{n=1}^N B_n&=:b+\frac{\eta}{\sqrt{N}}\,,\\
\end{aligned}
\]
where $a:=\mathbb{E}[A_n]$, $b:=\mathbb{E}[B_n]$, and $\xi\sim\mathcal{N}(0,\sigma_a^2),\ \eta\sim\mathcal{N}(0,\sigma_b^2)$ follow the normal distribution with $\sigma_a^2:=\mathrm{Var}[A_n]$, $\sigma_b^2:=\mathrm{Var}[B_n]$, we have the expression \eqref{h_nu_N} for $\Bar{h}_{\nu,N}$ written as
\begin{equation}\label{h_nv_N_random}
    \Bar{h}_{\nu,N}=\frac{a+\xi/\sqrt{N}}{b+\eta/\sqrt{N}}\,.
\end{equation}
To avoid the difficulty caused by the singularity with a probable zero value at the denominator  in \eqref{h_nv_N_random}, we introduce the auxiliary function
\begin{equation}\label{f_x_y_with_epsilon}
f(x,y)=\frac{a+x}{\epsilon + g(b+y-\epsilon)}\,, 
\end{equation}
with a parameter $\epsilon\in(0,b/2)$ and a regularized piecewise function $g$ defined by
\[  
g(z)=\left\{
\begin{aligned}
&\quad 0\,,\ z<-\epsilon\,,\\
&\frac{(z+\epsilon)^3}{6\epsilon^2}\,,\ -\epsilon<z<0\,,\\
&\frac{(\epsilon-z)^3}{6\epsilon^2}+z\,,\ 0<z<\epsilon\,,\\
&\quad z\,,\ z>\epsilon\,.\\
\end{aligned}
\right.
\]
We have $g(z)\ge 0$ for all $z\in\mathbb{R}$. Let $g_\epsilon(z):=\epsilon + g(z-\epsilon)$, the function $f(x,y)$ defined in \eqref{f_x_y_with_epsilon} can be written as
\[  f(x,y)=\frac{a+x}{g_\epsilon( b+y)}\,.\]
 A schematic plot of function $g_\epsilon(z)$ and a comparison between the functions $1/(b+y)$ and $1/g_\epsilon(b+y)$ is shown in Figure \ref{fig:plot_of_function_g}.

If we have $y+b-\epsilon>0$, then the two functions $1/(b+y)$ and $1/g_\epsilon(b+y)$ are identical. Given that the random variable $\eta\sim\mathcal{N}(0,\sigma_b^2)$ and that $\epsilon < b/2$, we have
\[  
\begin{aligned}
&\mathrm{Pr}\big( b+\frac{\eta}{\sqrt{N}}-\epsilon \le 0  \big) 
\\ = &\mathrm{Pr}\big( \eta \leq -\sqrt{N}(b-\epsilon) \big)\\
 = &\big\{ \eta\sim\mathcal{N}(0,\sigma_b^2) \Leftrightarrow \eta = \sigma_b\,\omega,\textup{ where }\omega \sim\mathcal{N}(0, 1) \big\}\\
 = &\mathrm{Pr}\big( \omega\leq-\frac{\sqrt{N}(b-\epsilon)}{\sigma_b} \big)\\
 = & \Phi\big( -\frac{\sqrt{N}(b-\epsilon)}{\sigma_b} \big)\,,
\end{aligned}
\]
where $\Phi$ is the cumulative distribution function of standard normal distribution. Recalling from Table \ref{Table_approx_Z_HO} for our numerical experiments with inverse temperature $\beta=2$, by employing a sample size $N=2^{28}$ to approximate the partition function $\Bar{Z}(\beta)$ we obtain a 95\% relative confidence interval $CI\approx 8\%$, which gives
\[ CI = 1.96\frac{\sigma_b}{\sqrt{N}b}\approx 0.08 \Longrightarrow \frac{\sqrt{N}b}{\sigma_b}\approx 24.5\,. \]
As the parameter $\epsilon$ satisfies $0<\epsilon<b/2$, we obtain
\[ \Phi(-24.5)\approx\Phi\big( -\frac{\sqrt{N}b}{\sigma_b} \big) < \Phi\big( -\frac{\sqrt{N}(b-\epsilon)}{\sigma_b} \big)< \Phi\big( -\frac{\sqrt{N}b}{2\,\sigma_b} \big)\approx\Phi(-12.3)\,, \]
which gives a tiny probability for a non-zero difference between $1/(b+\frac{\eta}{\sqrt{N}})$ and $1/g_\epsilon(b+\frac{\eta}{\sqrt{N}})$, since the probability $\mathrm{Pr}\big( \eta \leq -\sqrt{N}(b-\epsilon) \big)$ decreases fast as $N$ increases.

Further we obtain 
\[ 
\begin{aligned}
    & \frac{\partial f}{\partial x}=\frac{1}{\epsilon + g(b+y-\epsilon)},\quad 
      \frac{\partial f}{\partial y}=-\frac{(a+x)(g'(b+y-\epsilon))}{(\epsilon + g(b+y-\epsilon))^2}\,,\\
    & \frac{\partial^2 f}{\partial x^2}=0,\quad \frac{\partial^2 f}{\partial x\partial y}=-\frac{g'(b+y-\epsilon)}{(\epsilon + g(b+y-\epsilon))^2}\,,\\
    & \frac{\partial^2 f}{\partial y^2}=-\frac{(a+x)\big(g''(b+y-\epsilon)(\epsilon + g(b+y-\epsilon))+2\,g'(b+y-\epsilon)^2\big)}{(\epsilon + g(b+y-\epsilon))^3}\,. \\
\end{aligned}
\]
and the third derivatives of $f$ are bounded by $\mathcal{O}(\frac{1}{\epsilon^4})$. Applying Taylor's theorem we have
\[ 
\begin{aligned}
f(x,y)=&\, f(0,0)+f'_x(0,0)\,x+f'_y(0,0)\,y+f''_{xx}(0,0)\,\frac{x^2}{2}+f''_{yy}(0,0)\,\frac{y^2}{2}+f''_{xy}(0,0)\,xy\\
&+\mathcal{O}\big(\frac{1}{\epsilon^4}(|x|^3+|y|^3)\big)\,.  
\end{aligned}
\]
 With the condition $0<2\epsilon<b$ and substituting respectively $x$ and $y$ with $\xi/\sqrt{N}$ and $\eta/\sqrt{N}$, based on \eqref{h_nv_N_random} we further obtain, 
\begin{equation}\label{h_nu_bar_N_Taylor}
\Bar{h}_{\nu,N} \simeq \frac{a}{b}+\frac{1}{b}\frac{\xi}{\sqrt{N}}-\frac{a}{b^2}\frac{\eta}{\sqrt{N}}-\frac{\xi\eta}{b^2N}+\frac{a\eta^2}{b^3N}+\mathcal{O}\Big( \frac{1}{\epsilon^4}\big(\frac{|\xi|^3}{N^{\frac{3}{2}}}+\frac{|\eta|^3}{N^{\frac{3}{2}}}\big)\Big).
\end{equation}
Using the convexity of the functions $\xi^3$ and $\eta^3$, we apply Jensen's inequality to obtain
\[ \Big( \mathbb{E}\Big[ \frac{|\xi|^3}{N^{\frac{3}{2}}} \Big]\Big)^{\frac{4}{3}}\leq \mathbb{E}\Big[ \frac{|\xi|^4}{N^2} \Big]=\mathcal{O}\big( \frac{1}{N^2} \big)\,, \quad \Big( \mathbb{E}\Big[ \frac{|\eta|^3}{N^{\frac{3}{2}}} \Big]\Big)^{\frac{4}{3}}\leq \mathbb{E}\Big[ \frac{|\eta|^4}{N^2} \Big]=\mathcal{O}\big( \frac{1}{N^2} \big)\,,  \]
which yields
\[ \mathbb{E}\Big[ \frac{|\xi|^3+|\eta|^3}{N^{\frac{3}{2}}}\Big]=\mathcal{O}\big( \frac{1}{N^{\frac{3}{2}}} \big) \,. \]
Therefore, taking the expected value of $\Bar{h}_{\nu,N}$ based on \eqref{h_nu_bar_N_Taylor}, we obtain
\begin{equation}\label{E_h_bar_nu_N}
\mathbb{E}[\Bar{h}_{\nu,N}]=\frac{a}{b}-\mathbb{E}\big[ \frac{\xi\eta}{ab} \big]\frac{a}{b}\frac{1}{N} +\mathbb{E}\big[ \frac{\eta^2}{b^2} \big]\frac{a}{b}\frac{1}{N}+\mathcal{O}\big(\frac{1}{\epsilon^4 N^{\frac{3}{2}}} \big)\,. 
\end{equation}
Using specifically the covariance between $A_n$ and $B_n$,
\[ 
\begin{aligned}
\mathbb{E}[\xi^2]&=\mathbb{E}[(A_n-a)^2]=\sigma_a^2\,, \\
\mathbb{E}[\eta^2]&=\mathbb{E}[(B_n-b)^2]=\sigma_b^2\,, \\
\mathbb{E}[\xi\eta]&=\mathbb{E}[(A_n-a)(B_n-b)]=:\sigma_{ab}^2\,, \\
\end{aligned}
\]
we further obtain
\begin{equation}\label{Var_h_bar_nu_N}
\begin{aligned}
    \mathbb{E}[(h_\nu - \Bar{h}_{\nu,N})^2 ] &= \mathbb{E}\Big[ \Big( \frac{\xi}{a\sqrt{N}}\frac{a}{b}-\frac{\eta}{b\sqrt{N}}\frac{a}{b}+\mathcal{O}\big( \frac{1}{\epsilon^3}(\frac{\xi^2}{N} +\frac{\eta^2}{N})\big) \Big)^2 \Big]\\
    & = \big( \frac{a}{b} \big)^2\frac{1}{N}\big( \frac{\sigma_a^2}{a^2}+\frac{\sigma_b^2}{b^2}-2\frac{\sigma_{ab}^2}{ab} \big)+\mathcal{O}\big(\frac{1}{\epsilon^3 N^{\frac{3}{2}}} \big)\\
    & = \big\{ \,\textup{letting }\,\frac{\sigma_a^2}{a^2} + \frac{\sigma_b^2}{b^2}-\frac{2 \sigma_{ab}^2}{ab}=:\sigma_h^2 \big\}\\
    & = \frac{1}{N}\big( \frac{a}{b} \big)^2\sigma_h^2+\mathcal{O}\big(\frac{1}{\epsilon^3 N^{\frac{3}{2}}} \big).
\end{aligned}
\end{equation}

\begin{figure}
    \centering
    \includegraphics[width=0.9\textwidth]{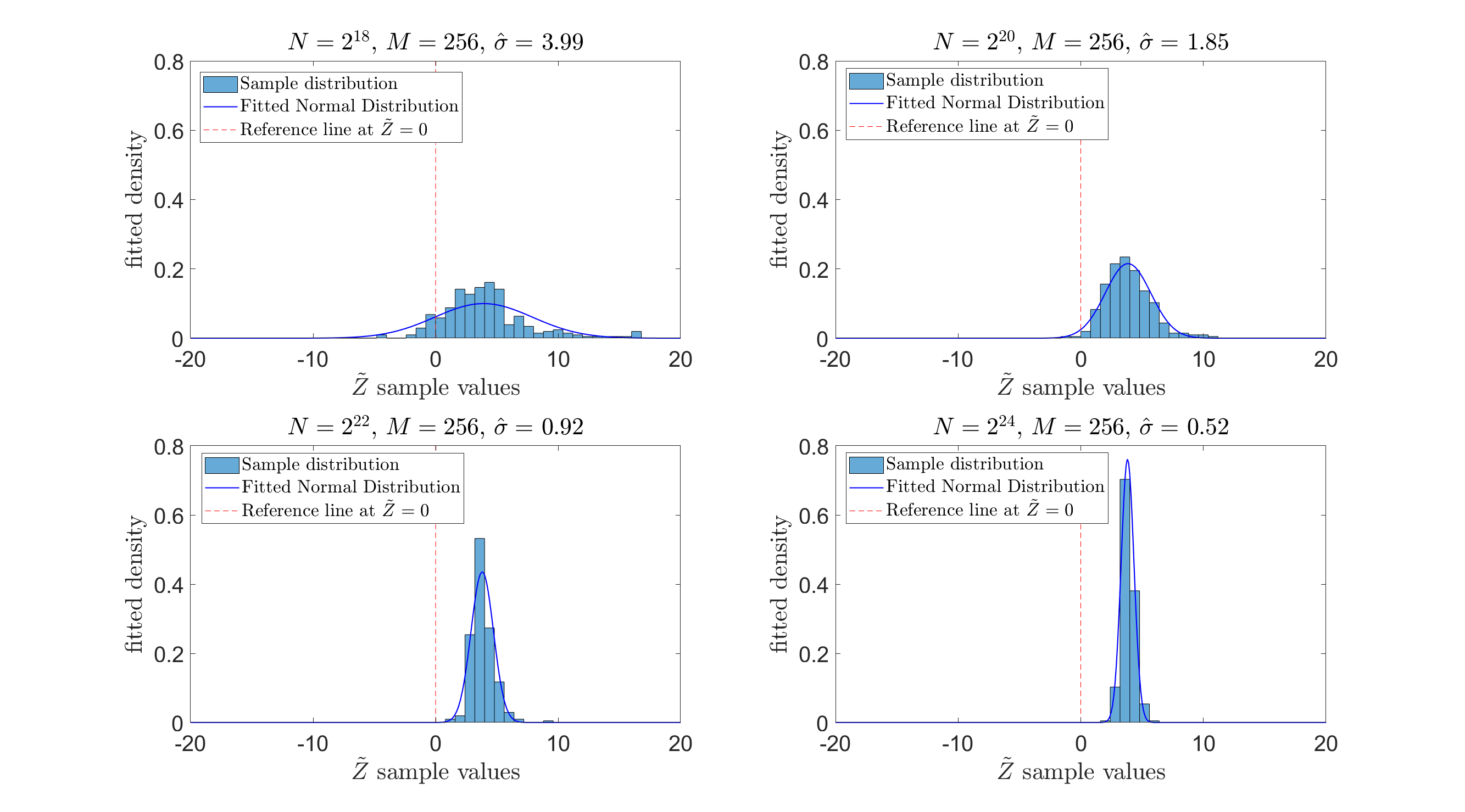}
    \caption{Histograms and the fitted normal distribution density functions for the partition function estimators $\Tilde{Z}_N$ with varying sample size $N=2^{18},2^{20},2^{22}$ and $2^{24}$ for the case \ref{Case_V1} under harmonic oscillator potential, with $n=6$, $d=3$, $\beta=2$, $\Delta t = 0.025$. For each fixed sample size, a total number of $M=256$ independent replicas of estimators $\{\Tilde{Z}_N^{(m)}\}_{m=1}^M$ are generated, and a sample standard deviation $\Hat{\sigma}$ is evaluated based on these replicas.}
    \label{fig:histogram_Z_N_estimator}
\end{figure}

According to the Berry-Esseen Theorem \cite{Probability_Berry_Esseen}, for independent identically distributed random variables $X_1$, $X_2$, ... with $\mathbb{E}[X_1]=0$, $E[X_1^2]=\sigma^2$, $\mathbb{E}[|X_1|^3]=\rho<+\infty$, let $Y_n=\frac{1}{n}\sum_{i=1}^n X_i$, then the cumulative distribution function (c.d.f.) $F_n$ of the normalized mean $\frac{\sqrt{n}Y_n}{\sigma}$ asymptotically approaches the c.d.f. of a standard normal distribution $\Phi$ as sample size $n$ increases. Specifically, $F_n$ satisfies
\[ \sup_{n,x}|F_n(x)-\Phi(x)|\leq \frac{C\rho}{\sigma^3\sqrt{n}}, \quad C<\frac{1}{2}\,.\]
Hence the stochastic variables $\xi$ and $\eta$ follow asymptotically normal distribution, with the difference between their c.d.f and the c.d.f. of their corresponding normal distributions bounded by $\mathcal{O}\big(\frac{1}{\sqrt{N}}\big)$, where $N$ is the sample size for our mean-field estimator $\Bar{h}_{\nu,N}$. Practically, we also test with a series of sample sizes $N=2^{18}$, $N=2^{20}$, $N=2^{22}$, and $N=2^{24}$. For each fixed sample size value, we generate $M=256$ replicas of independent estimators $\{\Tilde{Z}_N^{(m)}\}_{m=1}^M$, each obtained by multiplying the partition function estimator ${\Bar{Z}}_N^{(m)}$ with the constant factor $n!(2\pi\beta)^{\frac{3n}{2}}$. By fitting a normal distribution of the obtained estimators for different $N$ values, we indeed observe the empirical standard deviation $\Hat{\sigma}$ decreases with a factor around $1/2$ as the sample size $N$ increases with a factor $4$, as can be observed in Figure \ref{fig:histogram_Z_N_estimator}. 

Moreover, the two subfigures corresponding to $N=2^{22}$ and $N=2^{24}$ show that out of the $M=256$ independent replicas, no observations with a negative value of $\tilde{Z}_N$ are present, validating our previous assertion that $\mathrm{Pr}\Big(\big| \frac{1}{\big(b+\frac{\eta}{\sqrt{N}}\big)} - \frac{1}{g_\epsilon\big(b+\frac{\eta}{\sqrt{N}}\big)}\big|>0\Big)$ is negligible with a sufficiently large sample size $N$. Particularly in our numerical experiments, no bias related to the $\epsilon$-regularization \eqref{f_x_y_with_epsilon} is introduced.

\begin{figure}[h] 
\centering
\includegraphics[width=1\textwidth]
{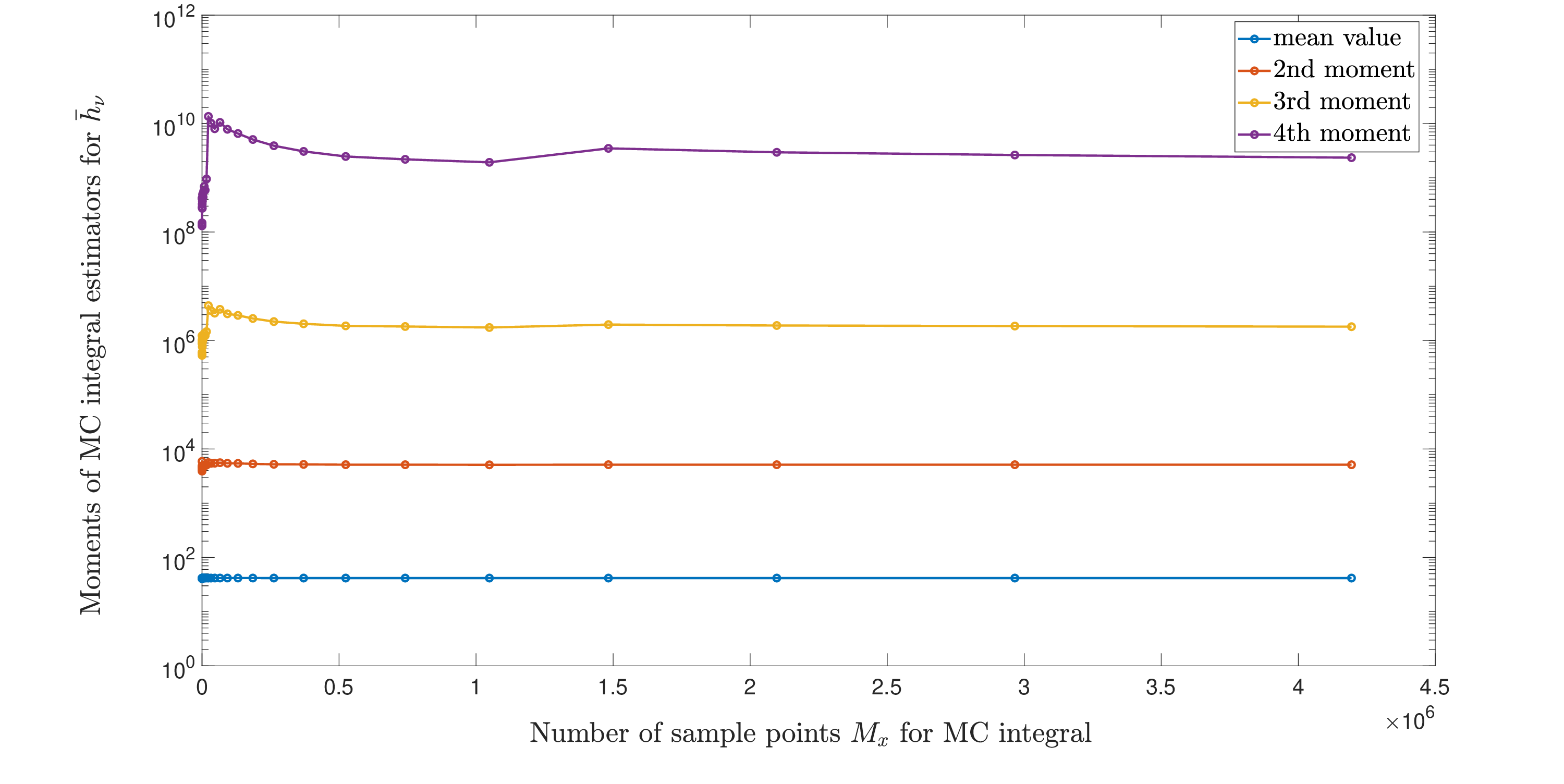}
\caption{ The mean and central moments up to order 4 of the mean-field energy estimators $\Bar{h}_\nu$ for $n=6$ fermions in dimension $d=3$, with $\beta=0.5$. The sample size $M_x=2^{22}$.}
\label{fig:Moments_of_h_nu}
\end{figure}

In order to check that the conditions for applying the Berry-Esseen theorem are satisfied and to verify that based on our sample set of size $N$, our estimator for the mean-field energy $\Bar{h}_{\nu,N}\simeq \frac{a}{b}+\frac{1}{\sqrt{N}}\big( \frac{1}{b}\xi-\frac{a}{b^2}\eta \big)$ is approximately normally distributed,  we also numerically evaluate the sample central moments of our estimators of $\mathcal{Z}_\nu$ and $\Bar{h}_\nu$. As can be observed from Figure \ref{fig:Moments_of_h_nu}, the moments gradually stabilize as the sample size $M_x$ increases, and remain bounded up to the $4$-th central moment. Combining the above reasoning with \eqref{E_h_bar_nu_N} and \eqref{Var_h_bar_nu_N}, we obtain
\begin{equation}
    \Bar{h}_{\nu,N}\simeq \frac{a}{b}\Big( 1\pm\frac{\sigma_h}{\sqrt{N}}+\mathcal{O}\big( \frac{1}{N}\big) \Big)\,.
\end{equation}
In practice, we use the sample mean and covariance values
\[ 
\begin{aligned}
    &\Bar{a}=\frac{1}{N}\sum_{n=1}^N A_n\,,\quad
    \Hat{\sigma}_a^2 = \frac{1}{N-1}\sum_{n=1}^N (A_n-\Bar{a})^2\,,\\
    &\Bar{b}=\frac{1}{N}\sum_{n=1}^N B_n\,,\quad
    \Hat{\sigma}_b^2 = \frac{1}{N-1}\sum_{n=1}^N (B_n-\Bar{b})^2\,,\\
    &\Hat{\sigma}_{ab}^2 = \frac{1}{N-1}\sum_{n=1}^N (A_n-\Bar{a})(B_n-\Bar{b})\,,\\
    &\Hat{\sigma}_h^2 = \frac{\Hat{\sigma}_a^2}{\Bar{a}^2} + \frac{\Hat{\sigma}_b^2}{\Bar{b}^2}-\frac{2\Hat{\sigma}_{ab}^2}{\Bar{a}\Bar{b}}\,,\\
\end{aligned}
\]
to obtain an empirical $95\%$ confidence interval for the mean-field estimator $\Bar{h}_{\nu,N}$ by 
\[\Big( \frac{\Bar{a}}{\Bar{b}}\big( 1 - 1.96\frac{\Hat{\sigma}_h}{\sqrt{N}}\big)\,,\ \frac{\Bar{a}}{\Bar{b}}\big( 1 + 1.96\frac{\Hat{\sigma}_h}{\sqrt{N}}\big)\Big).\]

To ensure the reliability of the above statistical method (Method 1), we also implement an alternative statistical method (Method 2). For simplicity, the specific ideas for the Method 2 can be explained with the following example:  if we want to employ a total sample size $M_x=2^{26}$, we first obtain $M_{1}=2^{10}$ independent estimators of the mean-field energy $\{\Bar{h}_\nu^{(i)}\}_{i=1}^{M_{1}}$, and each of these estimators $\Bar{h}_\nu^{(i)}$ is evaluated based on the independent sample sets $\{\mathbf{x}_0^{(i,j)},\mathbf{B}^{(i,j)}\}_{j=1}^{M_{2}}$, where $M_{2}=2^{16}$ so that the total sample size $M_x=M_1\times M_2$. By applying the Central Limit Theorem and fitting the obtained data $\{\Bar{h}_\nu^{(i)}\}_{i=1}^{M_{1}}$ with Normal distribution, we obtain the statistical confidence intervals for $\Bar{h}_\nu$, and they are in good consistency with the confidence intervals obtained with Method 1 for sufficiently large $M_1$ and $M_2$. 
Additional details for this implementation, along with sample code, are available on our open-access GitHub repository at \url{https://github.com/XinHuang2022/Path_Integral_MD_mean_field.git}.

\section{Particles with spin}\label{appendix_spin}
When spin is included a transposition $T_{ij}$ requires the coordinates for the  particle $i$, 
in particular $(x_i(t),y_i(t),z_i(t) ,s_i)$,
and particle $j$ %
to have the same spin, $s_i=s_j$. 
If $s_i\ne s_j$ thev particles $i$ and $j$ are treated as distinguishable, since the Hamiltonian $\widehat H$ does not depend on the spin. 
We obtain
two subsets of indistinguishable particles, one with spin $1/2$ and one with spin $-1/2$, and  determine the partition function for the system as follows. Let $\mathbf x=(\mathbf x_+,\mathbf x_-)$, where $\mathbf x_+=(x_1,y_1,z_1,s_1,\ldots, x_{n_+},y_{n_+},z_{n_+},s_{n_+})$ is the coordinate for particles with spin $s_1=s_2=\ldots=s_{n_+}=1/2$ and $\mathbf x_-=(x_{n_+ +1 },y_{n_+ +1 },z_{n_+ +1 },s_{n_+ +1 },\ldots, x_{n},y_{n},z_{n},s_{n})$ is the coordinate for particles with spin $-1/2$. An allowed coordinate permutation 
can then be written
$\mathbf x^\sigma=(\mathbf x_+^{\sigma^+},\mathbf x_-^{\sigma^-})$ 
where $\mathbf x_+^{\sigma^+}$ is the $\sigma^+$ permutation among the particles with spin $1/2$ and $\mathbf x_-^{\sigma^-}$ is the $\sigma^-$ permutation among the
particles with spin $-1/2$. The antisymmetric wave functions have the representation
\[\phi(\mathbf x)=
\sum_{\sigma^{+}\in\mathbf S^+}\sum_{\sigma^-\in\mathbf S^-}\frac{{\rm sgn}(\sigma^+)}{n_+!}\frac{{\rm sgn}(\sigma^-)}{n_-!}
\tilde\phi(\mathbf x^\sigma)\,,
\]
for $\tilde\phi\in L^2(\rset^{3n})$,
which yields the partition function
\[
\TR(e^{-\beta H_e})= \int_{\rset^{3n}}
\sum_{\sigma^{+}\in\mathbf S^+}\sum_{\sigma^-\in\mathbf S^-}\frac{{\rm sgn}(\sigma^+)}{n_+!}\frac{{\rm sgn}(\sigma^-)}{n_-!}
\mathbb E[e^{-\int_0^\beta V(\mathbf x_t,\mathbf X){\mathrm{d}}t}\delta(\mathbf x_\beta-\mathbf x_0^{(\sigma^+,\sigma^-)})]{\mathrm{d}}\mathbf x_0\,,
\] 
where
$n_\pm$ particles have spin $\pm1/2$ and the set of all coordinate permutations is denoted $\mathbf S^\pm$.
In the case of a separable potential 
\[
V(\mathbf x,\mathbf X) = \sum_{k=1}^{n_+} \tilde V(\mathbf x^+_{k},\mathbf X)
+\sum_{k=1+n_+}^{n} \tilde V(\mathbf x_{k}^-,\mathbf X)
\]
the partition function becomes
\[
\TR(e^{-\beta H_e})= \int_{\rset^{3n}}
=\int_{\rset^{3n}} \frac{{\rm det}\big(\mathcal W^+(\mathbf x^+_0)\big){\rm det}\big(\mathcal W^-(\mathbf x^-_0)\big)}{
n_+!n_-!(2\pi\beta)^{3n/2}} {\mathrm{d}}\mathbf x_0\,,
\]
where $\mathcal W^+$ is the $n^+\times n^+$ matrix with components 
\[
\mathcal W_{k\ell}^+\big(\mathbf x^+\big) := 
e^{-|\mathbf x_k^+(0)-\mathbf x_\ell^+(0)|^2/(2\beta)}
e^{-\int_0^\beta \tilde V(\mathbf B_k(t) + (1-\frac{t}{\beta}) \mathbf x_k^+(0)+ \frac{t}{\beta} \mathbf x_\ell^+(0),\mathbf X){\mathrm{d}}t}\,,
\]
and similarly for $\mathcal W_-$.

\end{document}